\documentclass{article}
\usepackage[T1]{fontenc}
\usepackage{lmodern}
\usepackage{mathtools, amsmath, amssymb, xcolor, amsthm, float, animate, mathrsfs, algorithm2e, csquotes, enumerate, dsfont}
\usepackage{booktabs}
\usepackage{pdfpages}
\usepackage[normalem]{ulem}
\mathtoolsset{showonlyrefs}
\usepackage{url}
\usepackage{subcaption}
\usepackage{hyperref}
\usepackage[export]{adjustbox}
\usepackage{geometry}
\usepackage{comment}
\usepackage{enumitem}

\newcounter{dummy} 
\numberwithin{dummy}{section}

\newtheorem{remark}[dummy]{Remark}
\newtheorem{theorem}[dummy]{Theorem}
\newtheorem{lemma}[dummy]{Lemma}
\newtheorem{example}[dummy]{Example}

\newtheorem{corollary}[dummy]{Corollary}
\newtheorem{proposition}[dummy]{Proposition}

\usepackage{appendix}

\DeclareMathOperator{\R}{\mathbb R}
\DeclareMathOperator{\N}{\mathbb N}

\DeclareMathOperator*{\argmin}{arg\,min}
\renewcommand{\d}{\mathop{}\!\mathrm{d}}

\usepackage[most]{tcolorbox}
\newtcolorbox{takeaway}{
    enhanced,
    breakable,
    colback=white,
    frame hidden,
    borderline west={2pt}{0pt}{black!65},
    boxrule=0pt,
    arc=0pt,
    left=8pt,
    right=4pt,
    top=4pt,
    bottom=4pt,
    before skip=8pt,
    after skip=8pt
}

\title{Finite speed Ornstein-Uhlenbeck-Kac processes and their applications in generative modeling
}
\author{Richard Duong\footnote{Institute of Mathematics,
	   Technische Universität Berlin,
	   Stra{\ss}e des 17.\ Juni 136, 
	   10623 Berlin, Germany,
	   {\ttfamily\{duong, auffenberg, steidl\}@math.tu-berlin.de},
      \url{https://tu.berlin/imageanalysis}.
	} \footnote{Corresponding author.}
\and Janis Auffenberg\footnotemark[1] 
\and Gabriele Steidl\footnotemark[1] 
 }
\date{\today}

\begin{document}
\maketitle

\begin{abstract}
    We study the Ornstein-Uhlenbeck-Kac (OU-K) process, where instead of a Brownian motion, a finite-speed Kac process defines the driving process. The OU-K process can be analyzed through a rescaled, time-inhomogeneous Kac process, characterized by a switching rate from the Kac part, and a convergence rate from the Ornstein-Uhlenbeck (OU) part. 
    The law of this inhomogeneous Kac process describes the solution of a wave equation, the Euler-Poisson-Darboux (EPD) equation. 
    For Kac switching rates being multiples of the OU convergence rate, we can explicitly derive this EPD solution via the sum of (the derivative of) two Legendre polynomials, restricted between two Dirac points representing the wave fronts. 
    Using the framework of piecewise-deterministic Markov processes (PDMP), we describe the conditional velocity fields of the OU-K process via our polynomial-law-formula. This allows us to apply the OU-K process in a generative flow model. 
    Since the OU-K velocity fields prove to be bounded by the finite Kac speed, the usual time truncations of diffusion models become unnecessary.
    We demonstrate advantages over diffusion and flow matching models in terms of generation sharpness on singular data sets, while achieving competitive results on CIFAR-10 image generation.
\end{abstract}

\section{Introduction}
\begin{figure}[ht] 
    \centering
    \includegraphics[width = 0.75\textwidth]{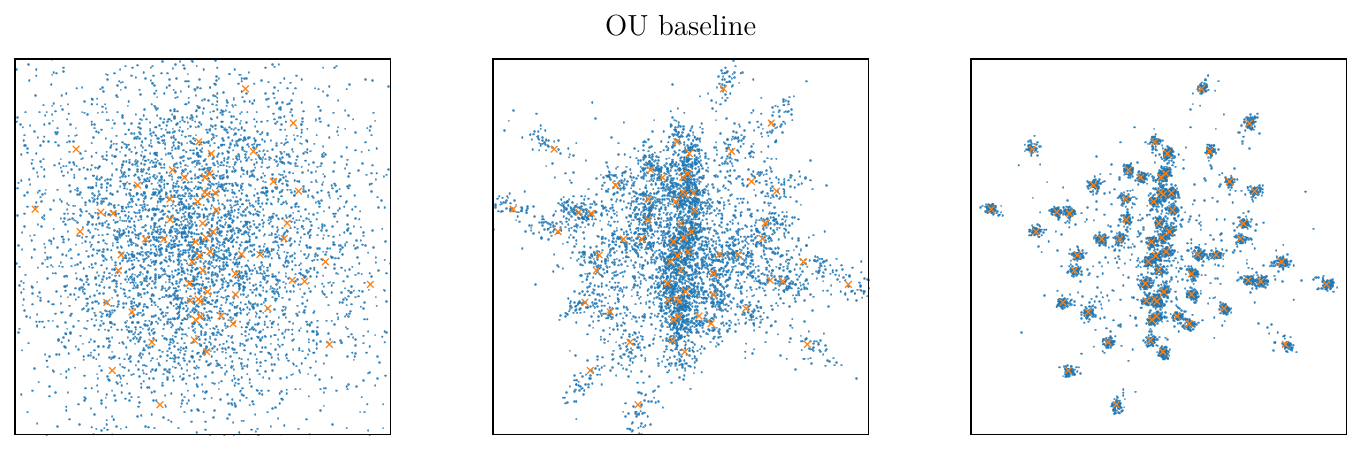}
    \\[2ex] 
    \includegraphics[width = 0.75\textwidth]{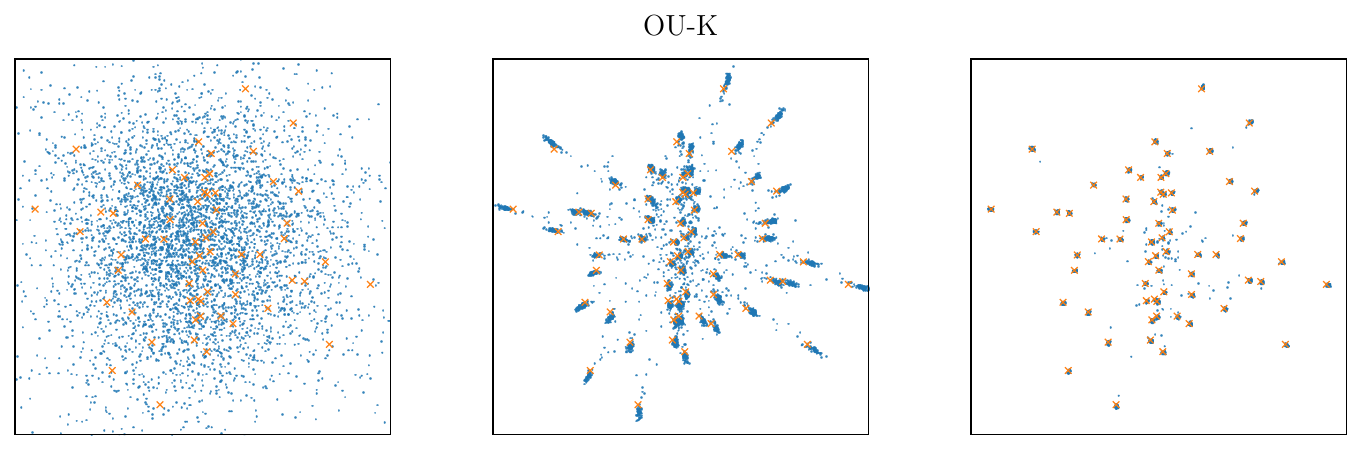}
       
    \caption{Generated samples after 15.000 (left), 30.000 (middle) and 100.000 (right) training steps of a polytope in dimension $d=32$ via the OU diffusion baseline (top) and our OU-K model (bottom), 
    orthogonally projected onto $\R^2$. By its better regularity, the OU-K model is able to approximate the vertices notably more precisely after significantly fewer training steps of the neural network.}
    \label{fig:OU-K_vs_OU_intro}
\end{figure}

Few approaches have shaped modern generative modeling as profoundly as flow matching and diffusion models \cite{ABV2023, lipman2023flow, Song2021}, which have now established as central tools in the field. 
Both approaches typically rely on stochastic interpolations between a target data variable and a reference (latent) Gaussian, and use neural networks to learn the flow transporting Gaussian noise to the data distribution.
Particularly relevant for the case of \emph{singular} data, the addition of Gaussian noise immediately extends the support of the interpolating distributions to the \emph{whole} Euclidean space.
Conversely, the learned flow must therefore collapse this full-dimensional support onto the possibly singular, lower-dimensional support of the data.
It has been known that this manifests in possible blow-ups of the velocity fields or score functions of diffusion models near the terminal time \cite{P2022, duong2026telegraphersgenerativemodelkac}, which may lead to ``fuzzy'' or ``blurry'' generations of singular data, see Figure \ref{fig:OU-K_vs_OU_intro}.

This issue is closely related to the \emph{infinite} propagation speed of the underlying Brownian motion and (mathematical) diffusion. Originated in physics \cite{C1958,V1958,MW1996}, it has been proposed to replace the heat-diffusion equation, which possibly features faster-than-light propagation, by the {relativistic heat equation}, also known as \emph{telegrapher's equation}. This hyperbolic PDE models the evolution of a {damped wave} propagating with \emph{finite} speed. As for the deep connection between diffusion equation and Brownian motion, the telegrapher's equation describes the law of a stochastic process, the \emph{Kac} or \emph{Telegraph} process. This process features a finite propagation speed, as well as a fixed switching rate $a>0$ corresponding to the damping of the wave.

In diffusion models \cite{Song2021}, the aforementioned stochastic interpolation usually takes the form of an Ornstein-Uhlenbeck process, driven by a Brownian motion. This work replaces the driving motion by the finite-speed Kac process $(K_t)_t$, resulting in the \emph{Ornstein-Uhlenbeck-Kac} (OU-K) process given by
\begin{equation}
Y_{t} = e^{- \theta t} Y_0 + \int_{0}^{t} e^{-\theta (t-s)} \d{K_s},
\end{equation}
where the exponential decay of the data $X_0$ is governed by a mean-reversion rate $\theta > 0$.
The analysis of the OU-K process and its law is generally cumbersome \cite{ratanov2021}. Yet, by the exponential prefactor, it can be linked to a \emph{time-inhomogeneous} Kac process $(Z_t)_t$ whose law solves a damped wave equation now with a time-\emph{dependent} damping coefficient. 

While obtaining the explicit solution for general time-inhomogeneous wave equations remains challenging \cite{MasoliverWeiss1994, AngelaniGarra2019}, the coefficients of the wave equation related to $(Z_t)_t$ are not arbitrary: they form a time-shifted version of the \emph{Euler-Poisson-Darboux} (EPD) equation \cite{bresters1973, DeGregorioOrsingher2020, GarraOrsingher2016},
\begin{equation}
    \partial_{t t} p + \frac{2 a}{\theta(1+t)} \partial_{t}p ~=~  \partial_{xx} p.
\end{equation}
For the cases $a = k \theta$, $k\in \N$, we apply the Riemann method to solve this PDE and acquire an \emph{elementary} expression for the law of $(Z_t)_t$ involving the two Legendre polynomials $L_{k}, L_{k-1}$ in Theorem \ref{thm:main-polynomial-law}. To our best knowledge, this representation is new to the literature. It further allows us to derive the explicit velocity fields associated with the probability flow of $(Z_t)_t$, providing the basis for our OU-K generative model. The use of other PDEs for generative modeling has been explored in \cite{GenPhys2023}.

Although the processes $(Y_t)_t$ and $(Z_t)_t$ are not Markovian themselves, 
they can be extended to become \emph{piecewise deterministic Markov processes} (PDMP). By this formulation, we may apply techniques and results of the PDMP literature towards their analysis: 
we obtain the (non-Gaussian) long-time limit of the OU-K process which plays the role of the reference/latent distribution in diffusion models, and we generalize the connection of $(Z_t)_t$ to its wave equation to a distributional setting by probabilistic arguments. 
Moreover, we extend the well-known conditional-expectation formula for probability-flow velocity fields to the PDMP setting, requiring only piecewise regularity of the process, see Proposition \ref{MSP-velo}.

Combining the above results, we introduce the \emph{OU-K generative model}, which integrates the OU-K process as the backbone of a flow matching framework. 
We theoretically establish that, under the assumption of bounded data, the conditional OU-K velocity fields involved in the flow matching objective remain uniformly bounded by the finite Kac speed. Consequently, unlike diffusion models, whose loss functions may exhibit singularities near the initial and/or terminal times, the OU-K model \emph{intrinsically} avoids such blow-ups and does not require the conventional time truncation schemes used in diffusion models.

While this theoretical advantage may appear subtle, our numerical experiments reveal 
practical benefits of the OU-K model over classical diffusion and flow matching approaches. On high-dimensional \emph{singular} data sets, our quantitative evaluations demonstrate that the OU-K model generates \emph{significantly sharper} samples with less training time than the diffusion and flow matching baselines, see again Figure \ref{fig:OU-K_vs_OU_intro}.
Moreover, on CIFAR-10, the OU-K model achieves a competitive FID score of $3.48$, compared with $3.59$ and $3.50$ for the flow matching and OU baselines, respectively.

\paragraph{Outline of the paper.}
In Section \ref{sec:wasserstein} we give a basic introduction into the mathematical background of flow matching and diffusion models. Section \ref{subsec:OUK} formally introduces the OU-K process and its invariant distribution. In Secion \ref{sec:inhom-Kac} we explicitly obtain the law of a related inhomogeneous Kac process by means of the EPD equation. This enables the derivation of its velocity field in Section \ref{sbsec:vkac}, and the conditional OU-K velocity field in Section \ref{sbsec:vkac-cond}. Building on these findings, the numerical Section \ref{sec:generative-modeling} introduces the OU-K model and demonstrates its advantages in data generation tasks.
Finally, Appendix \ref{sec:appendix} collects the proofs of several results presented in the paper, while Appendix \ref{sec:appendix-B} contains implementation details and additional figures.

\section{Ornstein--Uhlenbeck  process and generative flows}\label{sec:wasserstein}

\textbf{The standard Ornstein--Uhlenbeck (OU) process} $(X_t)_t$ with mean-reversion rate $\theta > 0$, variance parameter $\sigma^2 > 0$ and mean $\mu \in \R^d$ is given by the SDE
   \begin{equation}\label{Wiener-SDE}
       \d X_t = \theta (\mu - X_t) \d t + \sigma \d W_t,
   \end{equation}
where $(W_t)_t$ is a standard Brownian motion. 
An application of It{\^o}'s product rule to $e^{\theta t}X_t$ gives the explicit representation
   \begin{align}\label{def:OU}
       X_t &= e^{-\theta t} X_0 + \mu(1 - e^{-\theta t}) + \sigma \int_0^t e^{-\theta(t-s)} \d W_s \nonumber \\
       &\overset{d}{=} 
       e^{-\theta t} X_0 + \mu(1 - e^{-\theta t}) + \sqrt{\frac{1 - e^{-2 \theta t}}{2\theta}} \, \sigma \, \Xi,
   \end{align}
where $\Xi \sim \mathcal{N}(0, I_d)$ is independent of $X_0$. Hence, the OU process converges in distribution in the long-time limit $t \to \infty$ to
\begin{equation}\label{OU-latent}
    X_\infty \sim \mathcal{N}\Big(\mu, \frac{\sigma^2}{2\theta} I_d\Big).
\end{equation}
Denoting by $p_t$  the density of $X_t$, the corresponding
Fokker--Planck equation reads as
\begin{align}\label{OU-FPE}
    \partial_t p_t
    = -\nabla_x \cdot \left(\theta(\mu-\cdot )p_t \right)
    + \frac{\sigma^2}{2}\Delta_x p_t
    =
    - \nabla_x\cdot\left(v_t\, p_t\right)
\end{align}
with the \emph{probability-flow velocity field}
\begin{equation}\label{OU-velocity}
    v_t
    \coloneqq \theta(\mu-\cdot)
    - \frac{\sigma^2}{2} \underbrace{\nabla_x\log p_t}_{s_t}
    \end{equation}
and \emph{score function}
$s_t \coloneqq \nabla_x \log p_t$.
Applying a smooth time schedule $g \colon [0,1) \to [0,\infty)$ of the form
\begin{equation}\label{schedule-integral}
    g(t) \coloneqq \int_0^t \beta(s) \d s, \qquad \text{with } \beta > 0, \quad g(0) = 0, \quad \lim_{t \to 1} g(t) = \infty,
\end{equation}
and setting $\theta = \frac{1}{2}$, $\mu = 0$ and $\sigma = 1$, 
the time-changed process $\tilde X_t \coloneqq X_{g(t)}$, $t \in [0,1)$, satisfies the (forward) SDE
\begin{equation}\label{forward-OU-SDE}
   \d \tilde X_t = - \frac{1}{2} \beta(t) \tilde X_t \d t + \sqrt{\beta(t)} \d W_t, \qquad \tilde X_0 \sim \rho_{\rm data}.
\end{equation}
By \cite{Anderson1982}, the time reversal $\tilde X_t^* \overset{d}{=} \tilde X_{1-t}$, $t \in (0,1]$, satisfies the reverse SDE
\begin{equation}\label{reverse-SDE}
   \d \tilde X_t^* = \beta(1 - t) \left(\frac{1}{2} \tilde X_t^* + s_{1 - t}(\tilde X_t^*) \right) \d t + \sqrt{\beta(1 - t)} \d W_t^*, \qquad \tilde X_0^* \sim \mathcal{N}\Big(0, I_d\Big),
\end{equation}
where $(W_t^*)_t$ is some Brownian motion, and the score function determines the backwards drift.
\\[1ex]
\textbf{Score-based diffusion models} approximate the score function using a neural network by conditioning on target samples. Then, new samples of the target distribution $\rho_{\rm data}$ are obtained by numerically solving the reverse SDE \eqref{reverse-SDE} starting in samples from the latent distribution \eqref{OU-latent} given as $\mathcal{N}(0, I_d)$. 
However, the OU diffusion process gives rise to a score function that may become \emph{unbounded} as the process approaches the data distribution, see \cite{P2022, duong2026telegraphersgenerativemodelkac}, potentially leading to numerical instabilities during the sampling process.
\\[1ex]
\textbf{Flow matching models} on the other hand aim to approximate the vector field \eqref{OU-velocity}, solving the equation \eqref{OU-FPE}. Generally, we address such vector fields 
solving the \emph{continuity equation}
\begin{align}\label{eq:ce}
\partial_t p_t + 
\nabla_x \cdot (p_t v_t)=0
\end{align}
with some absolutely continuous Wasserstein curve $(p_t)_t$ as \emph{velocity fields associated with} $(p_t)_t$, see \cite{BookAmGiSa05}.
Notice that the continuity equation \eqref{eq:ce} remains meaningful in the sense of distributions. Hence, it also allows curves of probability measures $p_t$ featuring singularities, for which the evaluation of the score $\nabla_x \log p_t$ is not possible. Importantly, this will be the case in our OU-K model.

Given suitable regularity conditions on the velocity field $v_t$, it is then known \cite[Proposition 8.1.8]{BookAmGiSa05} that
the \emph{flow ODE} 
\begin{align}\label{eq:flow_ode} 
\partial_t \varphi(t,x)=v_t(\varphi(t,x)), \quad \varphi(0,x)=x,
\end{align}
has a unique solution $\varphi:I\times\R^d\to \R^d$ satisfying the push-forward relation 
$p_t = \varphi(t,\cdot)_\sharp p_0$.
Therefore, once such velocity field $v_t$ is known, we can transport samples from a given latent distribution $p_0$ to some target distribution $p_1$ by numerically solving the flow ODE \eqref{eq:flow_ode} and applying $p_1 = \varphi(1,\cdot)_\sharp p_0$. 
As an important remark, we will switch the flow direction in this work, similar to the diffusion approach \eqref{forward-OU-SDE}: the latent distribution will be given as $p_1$ from which we flow backwards to the target distribution $p_0$ simply by integrating $-v_{1-t}$ in \eqref{eq:flow_ode}, and applying $p_0 = \varphi(1,\cdot)_\sharp p_1$.

Since the velocity field $v_t$ is generally unknown, flow matching models approximate it by averaging over \emph{conditional velocity fields}, obtained by conditioning the flow on target samples. 
More precisely, they train a neural network $v_t^\theta$ to approximate 
$v_t$ by minimizing the \emph{Conditional Flow Matching} (CFM) loss given by
\begin{align} \label{cfm_loss}
     \mathcal{L}_{\mathrm{CFM}}(\theta) = \mathbb{E}_{t \sim \mathcal{U}(0,1), \, x_0 \sim p_0, \, x \sim p_t(\cdot| x_0)} \left[ \left\| v^\theta(t,x) - v(t,x| x_0) \right\|^2 \right], 
\end{align} 
see \cite{lipman2023flow, WS2025}. By conditioning the path $p_t$ on a given initial target sample $x_0 \sim p_0$, the {conditional velocity field} $v(t,x| x_0)$ becomes accessible.
In particular, the rescheduled OU process $(\tilde X_t)$ in \eqref{forward-OU-SDE}, conditioned on $\tilde X_0 = x_0$, possesses a conditional velocity field given by
\begin{align}\label{OU-cond-velocity}
v_t^{\rm OU}(x\mid x_0) =
\frac{\beta(t) e^{- \frac{g(t)}{2}}}{2 (1-e^{- g(t)})}
\left(e^{- \frac{g(t)}{2}} x - x_0\right),
\end{align}
see \cite[Theorem 3]{lipman2023flow}.
Unfortunately, it admits a blow-up as $t\to 0$, i.e.\ as the process approaches the target: 
intuitively, the velocity field must contract the entire Gaussian distribution to the Dirac mass $\delta_{x_0}$ as $t \to 0$, requiring increasingly large velocities.
As in the diffusion case, blow-ups of the velocity field may lead to numerical instabilities and, in some cases, even mode collapse, see \cite{CWDS2025}. 

This limitation motivates us to investigate an alternative stochastic-process backbone: the \emph{Ornstein–Uhlenbeck–Kac} process, which evolves with finite speed. Using this process, we obtain velocity fields that remain \emph{uniformly bounded} throughout the generative process.
Lastly, we note that the OU-K process is closely related to PDMPs, for which time-reversal approaches have also been developed, see \cite{Bertazzi2024}. In this work, however, we focus on a flow matching approach, which provides a natural framework for constructing generative probability flows.

\paragraph{Definitions and notational conventions.} 
We denote the Wasserstein-$p$ space of probability measures on $\R^d$ with finite $p$-th moment as $\mathcal{P}_p(\R^d)$, and equip it with the Wasserstein-$p$ distance $W_p$.
For functions depending on time and space, we frequently use the notation $f_t(x) \coloneqq f(t,x)$. For probability curves $(p_t)_t$ induced by a stochastic process $(X_t)_t$ via $p_t = {\rm Law}(X_t)$, we also speak of velocity fields \emph{associated with the process} $(X_t)_t$. For probability measures $p$ possibly having a singular part, we understand $p(x)$ in the sense of distributions.

\section{Finite speed Ornstein--Uhlenbeck process}\label{sec:finite-OU-K}

In this section, we introduce our main objects of interest: the Kac process and its Ornstein-Uhlenbeck variant. 
In contrast to the Wiener process, the Kac process, also known as ``persistent random walk'', describes a particle moving at \emph{finite} speed. 
Its close connection to the \emph{telegrapher's equation} (or damped wave equation) was established in the foundational works of Goldstein \cite{Goldstein1951} and Kac \cite{KAC1974}.
Unlike the diffusion equation associated with Brownian motion, the telegrapher's equation respects the physical constraint of finite propagation speed, which has made the Goldstein-Kac model particularly relevant in physics \cite{Broadbridge2019, MW1996, RGHS2026, Weiss2002}.
More recently, the Kac process was employed as the backbone of a generative model \cite{duong2026telegraphersgenerativemodelkac} and a subsequent distillation model \cite{distillkac2025}.
First results on the Ornstein-Uhlenbeck variant of the Kac process have been obtained in \cite{ratanov2021, Ratanov2022}. But, to the best of our knowledge, it has not previously been explored in the context of generative modeling.

In Subsection \ref{subsec:OUK}, we introduce the  Ornstein-Uhlenbeck-Kac process. Since in contrast to \eqref{def:OU}, 
its time marginals are not directly accessible, we investigate a related inhomogeneous Kac process in Subsection \ref{sec:inhom-Kac}. By extending this process to a PDMP, we derive its Fokker-Planck equation which takes the form of a specifically damped wave equation, the Euler-Poisson-Darboux equation. Solving this PDE explicitly reveals the inhomogeneous Kac law in Theorem \ref{thm:main-polynomial-law}. 
Based on this formula, we obtain the corresponding velocity fields in Section \ref{sec:velocity-kac}, which are essential ingredients of the OU-K generative model introduced in Section \ref{sec:generative-modeling}. 

All constructions in Section \ref{sec:finite-OU-K} are developed in the one-dimensional setting. Later, in 
Section \ref{subsec:multi-d}, they will be canonically extended to higher-dimensions by applying them componentwise, making them a suitable framework for generative modeling, for further details see \cite{ CKDS2026, duong2026telegraphersgenerativemodelkac}.

\subsection{Ornstein--Uhlenbeck--Kac process}\label{subsec:OUK}

Fix positive constants $a, c > 0$, where $a$ has the role of a \emph{switching intensity} (or damping), and $c$ is the \emph{moving speed} of the process. Let $(N_t)_t$ be a Poisson (counting) process with parameter $a$, i.e.\ a random process with values $N_t \in \N_0$, starting in $N_0 \equiv 0$, with independent increments which satisfy 
$$
N_t - N_s \sim {\rm Poi}(a(t-s)).
$$
Define the \emph{symmetric} and \emph{(time-)homogeneous Kac process} $(K_t)_t$ by
\begin{equation}\label{hom-Kac-process}
    K_t = B_{\frac{1}{2}} c \int_0^t (-1)^{N_s} \d s,
\end{equation}
where $B_{\frac{1}{2}} \sim {\rm Unif}(\{-1, 1\})$ is a Rademacher distributed coin-flip.
The \emph{Ornstein-Uhlenbeck-Kac} (OU-K) process with \emph{mean} $\mu \in \R$ and (mean-reversion) \emph{rate} $\theta > 0$ is defined via the SDE
\begin{equation}\label{OU-K-dynamics}
    \d Y_t = \theta(\mu - Y_t) \d t + \d K_t.
\end{equation}
Given an initial random variable $Y_0$, it has the explicit representation
\begin{equation}\label{def:OU-K}
Y_{t} = e^{- \theta t} Y_0 + \mu (1 - e^{- \theta t}) + 
\int_{0}^{t} e^{-\theta (t-s)} \d K_s,
\end{equation}
where  $\d K_s = B_{\frac{1}{2}} c (-1)^{N_s} \d s$. 

It is well-known that the usual Kac process converges to a Wiener process under the Kac scaling \eqref{eq:kac-scaling}, see \cite{KAC1974}. The following proposition establishes this also for the OU variants.
Its proof in Appendix \ref{appendix-proofs-A1} combines the known Kac scaling $(K_t^{(n)})_{t \in [0, T]} \overset{d}{\longrightarrow} (\sigma W_t)_{t \in [0, T]}$ with the continuous mapping theorem.

\begin{proposition}[Relation between OU and OU-K]\label{prop:kac-scaling}
    Let $\theta, \sigma > 0$ and $(a_n)_{n \in \N}, (c_n)_{n \in \N}$ with $a_n, c_n >0$ satisfy     
    \begin{equation}\label{eq:kac-scaling}
        a_n \to \infty,\qquad  c_n \to \infty, \qquad \frac{c_n^2}{a_n} \to \sigma^2 \qquad \text{as }  n \to \infty.
    \end{equation}
   For every $n \in \N$, consider the OU-K process $(Y_t^{(n)})_{t \geq 0}$ with mean-reversion rate $\theta$, mean $\mu \in \R$, speed $c_n$ and switching intensity $a_n$ starting in $Y_0^{(n)} \sim p_0$ for all $n \in \mathbb N$.
    Then, for any $T > 0$, it holds 
    \begin{equation}
        (Y_t^{(n)})_{t \in [0, T]} \overset{d}{\longrightarrow} (X_t)_{t \in [0, T]} \quad \text{as } n \to \infty
    \end{equation}
    in $C([0, T])$, where $(X_t)_{t \geq 0}$ is the Ornstein-Uhlenbeck process \eqref{Wiener-SDE} starting in $X_0 \sim p_0$.
\end{proposition}

The following theorem describes the law of the long-time limit of the OU-K process.  The invariant distribution \eqref{invariant-OU-K} will serve as the latent distribution in our flow matching approach. 
In contrast to  the stationary distribution $\mathcal N(\mu, \frac{\sigma^2}{2\theta} I_d)$ of the usual OU process, 
the OU-K limiting distribution is \emph{not} Gaussian.

\begin{theorem}[Limiting distribution of OU-K]\label{thm:invariantkacou}
    Let $(Y_t)_t$ be an OU-K process \eqref{def:OU-K}. 
Then the following holds true:
\begin{itemize}
\item[\rm{i)}]
  The OU-K process admits a unique 
    invariant distribution  given by the density
    \begin{equation}\label{invariant-OU-K}
        p_{\infty}(y) = \frac{1}{2 \beta B(\alpha, \alpha)}\left(\frac{\beta + (y - \mu)}{2 \beta}\right)^{\alpha - 1}\left(\frac{\beta - (y - \mu)}{2 \beta}\right)^{\alpha - 1} \mathds{1}_{[\mu - \beta,  \mu + \beta]}(y),
    \end{equation}    
    where $B(\alpha, \gamma) \coloneqq \int_{0}^{1}t^{\alpha - 1}(1 - t)^{\gamma - 1} \d{t}$ is the Euler-Beta function and
    $\alpha \coloneq \frac{a}{\theta}$ and $\beta \coloneq \frac{c}{\theta}$.
    Its variance is  $\sigma_\infty^2 = \frac{c^2}{\theta(\theta + 2a)}$.
    \item[\rm{ii)}]
    Assume that $\mathbb{E}[|Y_0|^q] < \infty$ for some $q \ge 1$. Denoting $p_t = {\rm Law}(Y_t)$, we have  
    \begin{equation}\label{Wasserstein-Bound}
        W_q(p_t, p_\infty) \leq \mathbb{E}[|Y_0 - Y_\infty|^q]^\frac{1}{q} e^{- \theta t}  \quad \text{for all } t \geq 0.
    \end{equation}
    In particular, it holds 
    $Y_t \overset{d}{\to} Y_\infty$ for $t \to \infty$.
    \end{itemize}
\end{theorem}

\begin{proof}
Part i) is due to \cite[Theorem 1]{Ratanov2022}.
For part ii), we note that
exponential convergence results hold more generally for certain subclasses of piecewise deterministic Markov processes, see \cite{BenaimLeBorgneMalrieuZitt2012, CloezHairer2015}. However, we give a short derivation of the explicit bound \eqref{Wasserstein-Bound}. Consider the OU-K processes $(Y_t)_t$ and $(\tilde Y_t)_t$ driven by the \emph{same} counting process $(N_t)_t$, where (by abuse of notation) $Y_0 \coloneqq Y_0$ and $\tilde Y_0 \coloneqq Y_\infty$, i.e.\ $(\tilde Y_t)_t$ is initialized at the stationary distribution $p_\infty$. Note that since both processes are synchronized by the same process $(N_t)_t$ per choice, it follows
\begin{equation}
\mathbb{E}[|Y_t - \tilde Y_t|^q] = e^{- q\theta t}\mathbb{E}[|Y_0 - \tilde Y_0|^q].
\end{equation}
By stationarity of $p_\infty$, the joint distribution $P_{Y_t, \tilde Y_t}$ is a coupling between $p_t$ and $p_\infty$, which yields
\begin{align}
    W_q(p_t, p_\infty) 
    \leq \left(\int_{\R \times \R} |x - y|^q \d P_{Y_t, \tilde Y_t}(x, y)\right)^\frac{1}{q} 
    = \mathbb{E}[|Y_t - \tilde Y_t|^q]^\frac{1}{q} = e^{-\theta t} \mathbb{E}[|Y_0 - \tilde Y_0|^q]^\frac{1}{q}.
\end{align}
\end{proof}

Next, we establish a trapping property of the OU-K process that is of independent interest: 
almost surely, the OU-K path will be trapped in a fixed compact interval after a finite (random) time. This result crucially relies on the finite-speed property of the Kac process, and does not hold for the usual OU process.

\begin{proposition}[Trapping property of OU-K process]
    There exists an a.s. finite random time $T \geq 0$ such that $Y_{t} \in [\mu - \frac{c}{\theta}, \mu + \frac{c}{\theta}]$ almost surely for all $t \geq T$.
\end{proposition}

\begin{proof}
For the existence of the finite hitting time $T$ of the interval $[\mu - \frac{c}{\theta}, \mu + \frac{c}{\theta}]$, see \cite[Theorem 3.1]{ratanov2021}. Note that the process stays trapped in the interval for all $t \ge T$ by the dynamics \eqref{OU-K-dynamics}.
\end{proof}

The following remark will be addressed later in Proposition \ref{Z-Kac-is-MSP}.

\begin{remark}\label{remark:PDMP}
    In contrast to OU process, both the Kac process $(K_t)_t$ and OU-K process $(Y_t)_t$ are not Markov processes themselves. By extending them by their speed/direction component 
    $$I_t = \frac{1}{2} B_{\frac{1}{2}} (-1)^{N_t} + \frac{3}{2} ~ \sim \, {\rm Unif}(\{1, 2\}),$$ the joint process $(K_t, I_t)_t$ becomes a Markov process {\rm \cite[Chapter 12.1]{EthierKurtz1986}}, and so does $(Y_t, I_t)_t$. Further, both extended processes are examples of piecewise deterministic Markov processes.
\end{remark}

\subsection{OU--K and  
Euler--Poisson--Darboux  equation}\label{sec:inhom-Kac}

The time-marginals of the homogeneous Kac process admit a representation via Bessel functions, see \cite{MW1996, duong2026telegraphersgenerativemodelkac}.
However, although the OU-K's long-time limit can be explicitly described by \eqref{invariant-OU-K}, the time-marginals of the OU-K process appear to elude a closed-form representation, see \cite[Section 5]{ratanov2021}.
To circumvent this issue, we rewrite the OU-K process via a \emph{time-inhomogenous} Kac process starting in $0$, and derive its law  explicitly  in Theorem \ref{thm:main-polynomial-law} via the shifted EPD equation. 
This explicit expression crucially enables the derivation of the (conditional) OU-K velocity field needed in our flow matching approach.

For this, we consider time-dependent switching intensities $a(t) > 0$. 
A \emph{time-inhomogeneous} Poisson (counting) process $(\tilde N_t)_t$ with rate $a(t)$ is then defined as 
a random process $\tilde N_t \in \N_0$ starting in $\tilde N_0 \equiv 0$, having independent increments which satisfy 
$$\tilde N_t - \tilde N_s \sim {\rm Poi} \left(\int_s^t a(u) \d u \right).$$
Analogously as in \eqref{hom-Kac-process}, we define the \emph{symmetric} and \emph{(time-)inhomogeneous Kac process} $(Z_t)_t$ by
\begin{equation}\label{inhom-Kac-process-Zt}
    Z_t = B_{\frac{1}{2}} c \int_0^t (-1)^{\tilde N_s} \d s.
\end{equation}

A \emph{key insight} is that the OU-K process $(Y_t)_t$ in \eqref{def:OU-K} admits a representation via a time-transformed inhomogeneous Kac process $(Z_u)_u$ yielding a specific time-dependent rate $a(u)$.

\begin{proposition}[Representation of OU-K via inhomogeneous Kac process]\label{prop:transformation}
Let $(Y_t)_t$ be an OU-K process \eqref{def:OU-K}. Then it holds
\begin{equation}\label{relation-to-Z_t}
    Y_{t} = e^{- \theta t} Y_0 + \mu (1 - e^{- \theta t}) + \frac{e^{-\theta t}}{\theta} \, c \,  Z_{e^{\theta t} - 1},
\end{equation}
where $(Z_u)_{u \geq 0}$ is a symmetric, inhomogeneous Kac process with speed $c_Z = 1$ 
and switching intensity 
\begin{equation}\label{specific-switching-OU-K}
    a(u) = \frac{a}{\theta (1 + u)}.
\end{equation}
\end{proposition}

\begin{proof}
    By substituting $u = e^{\theta s} - 1$ in \eqref{def:OU-K}, we obtain 
    \begin{align}
        Y_{t} &= 
        e^{- \theta t} Y_0 + \mu (1 - e^{- \theta t}) + B_{\frac{1}{2}} c e^{- \theta t} \int_{0}^{t} e^{\theta s} (-1)^{N_s} \d{s}\\
        &= e^{- \theta t} Y_0 + \mu (1 - e^{- \theta t}) + B_{\frac{1}{2}} c \frac{e^{- \theta t}}{\theta} \int_{0}^{e^{\theta t}-1} (-1)^{\tilde{N}_{u}} \d{u},
    \end{align}
    with the time-rescaled Poisson counting process
   $
        \tilde{N}_{u} = N_{\frac{1}{\theta }\log(1 + u)}.
   $
    Since $(N_{t})$ is a homogeneous Poisson counting process with intensity $a > 0$, it holds for any $0 \leq s < t$  that
    \begin{align*}
        \tilde{N}_{t} - \tilde{N}_{s} &= N_{\frac{1}{\theta }\log(1 + t)} - N_{\frac{1}{\theta}\log(1 + s)} \\ 
                                                       &\sim \text{Poi}\left(\frac{a}{\theta}(\log(1 + t) - \log(1 + s))\right) \\ 
                                                       &= \text{Poi}\left(\frac{a}{\theta}\int_{s}^{t}\frac{1}{1 + u} \d{u}\right). 
    \end{align*}
    Hence, $(\tilde{N}_{u})$ is a time-inhomogeneous Poisson process with rate 
    $a(u) \coloneq \frac{a}{\theta(1 + u)}$, and the above calculation shows \eqref{relation-to-Z_t}.
\end{proof}

The representation \eqref{relation-to-Z_t} can be considered the analogue to the OU formula \eqref{def:OU}. 

As in the homogeneous case, see Remark \ref{remark:PDMP}, the 
inhomogeneous Kac process $(Z_t)_t$ gives rise to a  piecewise-deterministic Markov process (PDMP). In the following we want to formulate this result without going extensively into detail. We point to the rich PDMP literature, e.g.\ \cite{AzaisBardetGenadotKrellZitt2014, davis1984, davis1993}, for further references.

\paragraph{Piecewise deterministic Markov processes.}

Next, we discuss a special subclass of PDMPs, which we call \emph{Markov switching processes}. Such process  $(X, I)$ consists of two components, the \emph{continuous} or \emph{spatial} component $X=(X_t)_t$ and the \emph{discrete} component $I=(I_t)_t$, such that $(X_t, I_t) \in \R^d \times \{1,\dots,n\}$ for all $t \ge 0$.

The continuous component $X$ moves deterministically between random times. This deterministic movement is modeled by vector fields 
\begin{equation}
    F^{i}: \R^{d} \to \R^{d}, \quad i=1,\ldots,n,
\end{equation}
where we restrict ourselves to continuously differentiable fields $F^i$ satisfying the Lipschitz condition
    \begin{equation}\label{ass:Lipschitz}
        \|{F^{i}(x) - F^{i}(y)}\|_{2} \leq L^{i} \|{x - y}\|_{2} \quad \text{for all } x, y \in \R^{d}
    \end{equation}
    for some constants $L^{i} > 0$.
Assumption \eqref{ass:Lipschitz}  guarantees global existence of the solution $\phi^i$ of the corresponding Cauchy problem 
\begin{equation}\label{Cauchy-problem-phi}
    \frac{\d{}}{\d{t}}\phi_{t}^{i}(x) = F^{i}(\phi_{t}^{i}(x)), \qquad \phi_{0}^{i}(x) = x,
\end{equation}
satisfying $\phi_{t + s}^{i}(x) = \phi^{i}_{t}(\phi^{i}_{s}(x))$ and $(\phi_{s}^{i})^{-1}(x) = \phi_{-s}^{i}(x)$.

The random times between the deterministic movements of $X$ are governed by the \emph{switching rates}: for  $i, j \in \{1, \ldots, n\}$, we consider continuous, bounded,
time-dependent rates 
$\lambda^{ij}:\R_{\geq 0} \times \R^{d} \to \R_{\geq 0}$
with 
\begin{align}\label{ass:continuous-switch}
   \lambda^i \coloneqq \sum_{j=1}^n \lambda^{ij}  > 0, \quad \text{and} \quad \lambda^{ii}=0.
\end{align}
Here $\lambda^{ij}$ describes the rate at which the process changes its dynamics from mode $i$ to $j$, and $\lambda_i$ is the total outgoing rate from mode $i$.
To indicate the time-dependence, we also write $\lambda_t^{ij}(x) \coloneqq \lambda^{ij}(t,x)$.

A \emph{Markov switching process} $(X, I)$ on $\R^{d} \times \{1, \ldots n\}$ with characteristic tuple  $(F, \lambda)$ is a stochastic process, where the continuous part $X = (X_t)_t$ is determined by the dynamics  
\begin{equation}\label{eq:mspdef}
    \d{X_{t}} = F^{I_{t}}(X_{t})\d{t}, \quad t > 0,
\end{equation}
and the discrete part $I = (I_t)_t$ is given as follows: 
let $T_0=0$ and $T_k$, $k=1,2,\dots$, be the sequence of random switching times on $[0,\infty)$, where the conditional laws ${\rm Law}(T_{k+1} \mid X_{T_k}, I_{T_k}, T_k)$ are inductively defined by
\begin{equation}\label{switching-times}
    \mathbb{P}(T_{k + 1} > T_{k} + t \mid X_{T_k} = x, I_{T_k}=i) = \exp\left(-\int_{0}^{t} \lambda^{i}_{s+T_k}(\phi^{i}_{s}(x)) \d{s}\right).
\end{equation}
Here, $\phi_t^i(x)$ is given by the solution of the Cauchy problem \eqref{Cauchy-problem-phi}.
Then the law of the discrete part $I$ at the next switching time $T_{k+1}$ is given by
\begin{equation}
    \mathbb{P}(I_{T_{k + 1}} = j \mid X_{T_{k+1}} = \tilde x, I_{T_k}=i) = \frac{\lambda^{ij}_{T_{k+1}}(\tilde x)}{\lambda^{i}_{T_{k+1}}(\tilde x)}.
\end{equation}
We set $I_t \equiv I_{T_k}$ for all $t \in [T_k, T_{k+1})$. By the assumption \eqref{ass:continuous-switch}, it is guaranteed that $T_k \to \infty$ for $k \to \infty$ almost surely, such that $(X_t, I_t)$ is defined for all $t \in [0, \infty)$.

Every Markov switching process is indeed a PDMP, and hence a Markov process, see \cite{AzaisBardetGenadotKrellZitt2014} for the homogeneous case. This also translates to our time-inhomogeneous setting by the usual time-space augmentation.

The following proposition describes the OU-K and Kac process in terms of Markov switching processes. 

\begin{proposition}[Extension of Kac and OU-K process to PDMP]
\label{Z-Kac-is-MSP}
    Let $(Z_t)_t$ be a symmetric, inhomogeneous Kac process with speed $c_Z > 0$, continuous switching rate $a_Z(t) > 0$, and corresponding counting process $(N_t^Z)$.
    Further, let $(Y_t)_t$ be an OU-K process with speed $c_Y > 0$, switching intensity $a_Y > 0$, mean $\mu \in \R$ and rate $\theta > 0$, with counting process $(N_t^Y)$.\\
    Then, together with the discrete components
    \begin{equation}\label{switching-processes-from-poisson}
        I_t^Z = \frac{1}{2} B_{\frac{1}{2}} (-1)^{N_t^Z} + \frac{3}{2}, \quad ~
        I_t^Y = \frac{1}{2} B_{\frac{1}{2}} (-1)^{N_t^Y} + \frac{3}{2}
        \quad \sim ~ {\rm Unif}(\{1, 2\}),
    \end{equation}
    the pairs $(Z_t, I_t^Z)$ and $(Y_t, I_t^Y)$ define Markov switching processses with characteristic tuples $(F_Z, \lambda_Z)$ and $(F_Y, \lambda_Y)$, resp., given by
    \begin{align}\label{Z-Tuple}
        F_Z^i(x) &= (-1)^{i}c_Z, & F_Y^i(x) &= \theta(\mu - x) + (-1)^{i}c_Y, & i &= 1,2, \\
        \lambda_Z^{ij}(t,x) &= a_Z(t) \mathds 1_{\{j \ne i\}}, & \lambda_Y^{ij}(t,x) &= a_Y \mathds 1_{\{j \ne i\}}, & i, j &= 1, 2. 
    \end{align}
\end{proposition}

For convenience, the straightforward proof is given in Appendix \ref{appendix-proofs-A1}.

\paragraph{Fokker--Planck equation of Markov switching processes.}
For a Markov switching process $(X,I)$, we introduce the mode-wise measures
\begin{equation}\label{mode-wise-measure}
    \rho_{i,t} (A) \coloneqq \mathbb{P}(X_t \in A, \, I_t = i), \quad A \subset \R^d ~ \text{Borel-measurable},
\end{equation}
which are \emph{sub}probability measures, i.e.\ $0 \le \rho_{i,t}(\R^d) \le 1$ for all $i=1,\dots,n$, $t\ge 0$.
The mode-wise laws of a Markov switching process satisfy a Fokker-Planck type equation, the \emph{Kolmogorov forward equation}, where the diffusion term is replaced by a term describing the jumps of the discrete component.
For the proof, see Appendix \ref{appendix-proofs-A1}.

\begin{proposition}[Mode-wise KFE]\label{FPE-MSP}
    Let $(X,I)$ be a Markov switching process with characteristic tuple $(F,\lambda)$ satisfying \eqref{ass:Lipschitz} and \eqref{ass:continuous-switch}.
    Then, for every mode $i=1,\dots,n$, the mode-wise measures $\rho_{i,t}$ satisfy the Kolmogorov forward equation
    \begin{equation}\label{mode-wise-KFE}
        \partial_t \rho_{i,t} = - {\rm div} \big( F^i \, \rho_{i,t} \big)
        + \sum_{j \ne i} \left(\lambda^{ji}\rho_{j,t} - \lambda^{ij}\rho_{i,t} \right)
    \end{equation}
    in the sense of distributions on $(0, \infty) \times \R^d$. 
\end{proposition}

It is known that the law of Kac processes can be linked to the solution of a hyperbolic analogue of the Fokker-Planck equation, the \emph{damped wave equation}, see \cite{KAC1974} for the time-homogeneous case, and \cite{Kaplan1964} for the time-inhomogeneous case.
Note that both works assume smooth initial conditions. By viewing Kac processes as PDMPs, we can give a rather concise derivation of this result  for the general distributional setting.

\begin{theorem}[Inhomogeneous Kac process and damped wave equation]\label{prop:PDE-1}
Let $(Z_t)_{t \geq 0}$ be a symmetric, inhomogeneous Kac process with speed $c > 0$ and continuous switching rate $a(t) > 0$. 
Then, the law $p(t,x)$ of $(Z_t)_t$ satisfies the damped wave equation
\begin{align}\label{damped-wave}
    \begin{cases}
    \partial_{tt}p + 2 a(t) \partial_{t}p \,=\,  c^2 \partial_{xx} p,\\ 
    p(0, x) = \delta(x),\\ 
\partial_{t}p(0, x) = 0.\end{cases}
\end{align}
in the sense of distributions on $(0,\infty) \times \R$.
\end{theorem}

\begin{proof}
    By Proposition \ref{Z-Kac-is-MSP}, we know that $(Z_t, I_t^Z)$ defines a Markov switching process with the tuple $(F_Z, \lambda_Z)$ given in \eqref{Z-Tuple}.
    By symmetry of the switching rates $\lambda_Z^{ij}=\lambda_Z^{ji}=a(t)$ and $\partial_x ((-1)^{i}c) = 0$, Proposition \ref{FPE-MSP} yields the mode-wise Kolmogorov forward equations 
    \begin{equation}\label{KFE-mode}
        \partial_{t}\rho_{i, t} = (-1)^{i + 1}c \, \partial_{x}\rho_{i, t} + a(t) (\rho_{j, t} - \rho_{i, t}), \qquad i,j =1,2, \, j \neq i,
    \end{equation}
    for the mode-wise measures $\rho_{i, t}$ given in \eqref{mode-wise-measure}, in the sense of distributions.

    Note that the following calculations also hold in the sense of distributions, where differentiation acts on smooth, compactly supported test functions on $(0,\infty) \times \R$; for more details, see Appendix \ref{appendix-proofs-A1}.
    Using $p = \rho_{1, t} + \rho_{2, t}$ and the notation $j \coloneqq \rho_{2, t} - \rho_{1, t}$, we obtain
    \begin{align*}
        \partial_{t} p(t, x) &= \partial_{t}(\rho_{1, t} + \rho_{2, t}) 
                             = c \, \partial_{x} \rho_{1, t} + a(t) (\rho_{2, t} - \rho_{1, t}) - c \, \partial_{x} \rho_{2, t} + a(t)(\rho_{1, t} - \rho_{2, t}) \\ 
                             &= - c \, \partial_{x}(\rho_{2, t} - \rho_{1, t})  
                             = - c \, \partial_{x} j(t, x),\\
            \partial_{t} j(t, x) &= \partial_{t}(\rho_{2, t} - \rho_{1, t})  
                             = -c \, \partial_{x}\rho_{2, t} + a(t) (\rho_{1, t} - \rho_{2, t}) - c \, \partial_{x} \rho_{1, t} - a(t)(\rho_{2, t} - \rho_{1, t}) \\ 
                             &= - c \, \partial_{x} p(t, x) - 2 a(t) j(t,x).
    \end{align*}
    Combining both equations yields the damped wave equation
    \begin{align*}
        \partial_{tt}p(t, x) &= \partial_{t}(- c \, \partial_{x} j(t, x)) = -c \, \partial_{x}(\partial_{t}j(t, x)) = c^{2} \partial_{xx} p(t, x) + 2 a(t) c \, \partial_{x} j(t, x) \\ 
        &= c^{2} \partial_{xx}p(t, x) - 2 a(t) \partial_{t}p(t, x).
    \end{align*}
    The initial condition $p(0, \cdot) = \delta$ holds directly per definition \eqref{inhom-Kac-process-Zt}. Moreover, by symmetry it holds $\rho_{1, 0} = \rho_{2, 0} = \frac{1}{2} \delta$. Then, it follows 
    \begin{equation}
        \partial_{t} p(0, \cdot) = \partial_{t} \rho_{1, 0} + \partial_{t} \rho_{2, 0} = \frac{c}{2}\delta' - \frac{c}{2}\delta' = 0
    \end{equation}
    by evaluating \eqref{KFE-mode} at $t=0$.
\end{proof}

The treatment of the wave equation \eqref{damped-wave} in the distributional sense becomes necessary, since the inhomogeneous Kac law $p(t,x)$ admits singularities in form of Dirac measures at the wave fronts, as we will see below.

\paragraph{Explicit solution formula for the inhomogeneous Kac law.}

The advantage of using $(Z_u)_{u \geq 0}$ in Proposition \ref{prop:transformation} is that its probability law solves the damped wave equation \eqref{damped-wave}.
Yet, for general time-inhomogeneous Kac processes, the explicit representation of the solution of \eqref{damped-wave} remains challenging, and is known only for certain combinations of switching rate $a(t)$ and moving speed $c(t)$, see \cite{AngelaniGarra2019, MartinucciSpina2025, MasoliverWeiss1994}. 
Interestingly, for the specific switching intensity $a(u) = \frac{a}{\theta (1 + u)}$ in \eqref{specific-switching-OU-K}, the law is the solution of a \emph{weakly dissipative} wave equation, which can be analytically solved via the Fourier transform. 

\begin{proposition}[Law of inhomogeneous Kac process and its Fourier transform]\label{prop:PDE}
Let $(Z_u)_{u \geq 0}$ be a symmetric, inhomogeneous Kac process with speed $c_Z = 1$ and switching rate $a(u) = \frac{a}{\theta(1 + u)}$, where $a, \theta > 0$ are fixed. Then the following holds true.
\begin{itemize}
\item[\rm i)]    
The law of $Z_u$ admits a singular and absolutely continuous part via
\begin{equation}
    p_Z(u,x) = \frac{1}{2} \exp\left(-\int_0^u a(s) \d s \right)( \delta(x - u) + \delta(x + u))
    ~+~ p_Z^{\rm abs}(u,x),
\end{equation}
and $p_Z(u,x)$ satisfies the PDE
\begin{align}\label{dissipative-wave}
    \begin{cases}
    \partial_{uu}p + \frac{2 a}{\theta(1 + u)} \partial_{u}p \,=\,  \partial_{xx} p,\\ 
    p(0, x) = \delta(x),\\ 
\partial_{u}p(0, x) = 0.\end{cases}
\end{align}
\item[\rm ii)]   
Its Fourier transform\footnote{using the convention $\hat f(\xi) = \int_{-\infty}^\infty f(x) e^{-i \xi x} \d x$ and $f(x) = \frac{1}{2\pi} \int_{-\infty}^\infty \hat f(\xi) e^{i \xi x} \d \xi$.} is given as
\begin{equation}\label{FT-of-law}
    \hat p_Z(u, \xi) = -\frac{\pi}{2}  |\xi| (1+u)^\rho \Big(J_{\rho-1} (|\xi|) M_{\rho} ((1+u)|\xi|) - M_{\rho-1} (|\xi|) J_{\rho} ((1+u)|\xi|) \Big),
\end{equation}
where $\rho \coloneqq - \frac{1}{2} (\frac{2a}{\theta} -1)$, and $J_\rho, M_\rho$ are Bessel functions of first and second kind, respectively.
Furthermore, for non-integer $\rho$, we obtain a simplified expression for the Fourier transform
\begin{equation}\label{FT-of-law-simple}
    \hat p_Z(u, \xi) = \frac{\pi}{2 \sin(\rho \pi)}  |\xi| (1+u)^\rho \Big(J_{1-\rho} (|\xi|) J_{\rho} ((1+u)|\xi|) + J_{\rho-1} (|\xi|) J_{-\rho} ((1+u)|\xi|) \Big).
\end{equation}
\end{itemize}
\end{proposition}

\begin{proof}
    The result in i) follows from Theorem \ref{prop:PDE-1}.
    Statement ii) can be found in  \cite{Wirth2004}, where we corrected an erroneous sign in \eqref{FT-of-law-simple}.
\end{proof}


The inverse Fourier transform of \eqref{FT-of-law-simple}, understood in the distributional sense, formally reads as
\begin{equation}
    p_Z(u,x) = \frac{1}{2 \sin(\rho \pi)} (1+u)^\rho \int_0^\infty \xi
    \Big(J_{1 - \rho} (\xi) J_{\rho} \big((1+u)\xi \big) + J_{\rho-1} (\xi) J_{-\rho} \big((1+u)\xi \big) \Big)
    \cos(\xi x) \d \xi.
\end{equation}
Surprisingly, if the switching intensity $a$ is an \emph{integer multiple} of the rate $\theta$, the integral dissolves into elementary expressions, as the following examples demonstrate.

\begin{example}\label{explicit-examples-Fourier}
\begin{enumerate}
\item For $a=0 \cdot \theta$, i.e.\ $\rho = \frac{1}{2}$,
     we simply have two Dirac points 
    \begin{equation}
        p_Z(u, x) = \frac{1}{2} (\delta(x-u) + \delta(x+u)),
    \end{equation}
    with the Fourier transform
    \begin{equation}
       \hat p_Z(u, \xi) = \cos(u \xi).
    \end{equation}
\item For $a=1 \cdot \theta$, i.e.\ $\rho = -\frac{1}{2}$,
    we obtain a uniform distribution enclosed by two Dirac points, 
    \begin{equation}
        p_Z(u, x) = \frac{1}{2(1+u)} (\delta(x-u) +\delta(x+u)) + \frac{1}{2(1+u)} \mathds{1}_{\{|x| < u\}},
    \end{equation}
    with the Fourier transform
    \begin{equation}
       \hat p_Z(u, \xi) = \frac{1}{1+u} \Big( \cos(u \xi) + \frac{\sin(u \xi)}{\xi}\Big).
    \end{equation}  
\item For $a=2 \cdot \theta$, i.e.\ $\rho = -\frac{3}{2}$,
 we get a quadratic polynomial density enclosed by two Dirac points, 
\begin{equation}
        p_Z(u, x) = \frac{1}{2(1+u)^2}(\delta(x-u) + \delta(x+u)) + 
        \frac{3u^2 +6u + 4 -3x^2}{4(1+u)^3}  \mathds{1}_{\{|x| < u\}},
    \end{equation}
 with the Fourier transform
    \begin{equation}
       \hat p_Z(u, \xi) = \frac{1}{(1+u)^2} \Big( \cos(u \xi) + \frac{2+3u}{1+u} \cdot \frac{\sin(u \xi)}{\xi} + \frac{3}{1+u} \cdot \frac{\sin(u \xi) - u \xi \cos(u \xi)}{\xi^3}\Big).
    \end{equation}  
    \end{enumerate}
\end{example}

Driven by the above examples, we now treat the \emph{general} case of integer multiples $a=k\theta$, i.e.\ $\rho = \frac{1}{2} - k$, $k \in \N$. To this end, we consider the (shifted) \emph{Euler-Poisson-Darboux} (EPD) equation
\begin{align}\label{EPD}
    \begin{cases}
        \partial_{\tau \tau} p + \frac{2 k}{\tau} \partial_{\tau}p ~=~  \partial_{xx} p,\\ 
        p(1, x) = \delta(x),\\ 
        \partial_{\tau}p(1, x) = 0.\end{cases}
\end{align}
This yields an explicit representation of the law of $Z_u$ in terms of \emph{two Legendre polynomials}.
To the best of our knowledge, the representation given in \eqref{legendre-formula-abs} is new to the literature.

\begin{theorem}
[Law of inhomogeneous Kac process for $\frac{a}{\theta} \in \N$]
\label{thm:main-polynomial-law}
    Let $\frac{a}{\theta}=k \in \N$, and $(Z_u)_{u \ge 0}$ be a symmetric, inhomogeneous Kac process with speed $c_Z = 1$ and switching rate $a(u) = \frac{a}{\theta(1+u)} = \frac{k}{1+u}$.\\ 
    Then, the law of $Z_u$ reads as  
    \begin{equation}\label{legendre-formula-total}
        p_Z(u, x) = \frac{1}{2(1+u)^k}\big(\delta(x-u) + \delta(x+u) \big) + 
        P_Z^{(k)}(u,x) \mathds{1}_{\{|x| < u\}},
    \end{equation}
    where $P_Z^{(k)}(u, \cdot)$ is a polynomial of degree $2k-2$ given as
    \begin{equation}\label{legendre-formula-abs}
        P_Z^{(k)}(u,x) = \frac{1}{2(1+u)^k} \left( L_k'\Big(\frac{(1+u)^2 - x^2 + 1}{2(1+u)} \Big) - \frac{1}{1+u} L_{k-1}'\Big(\frac{(1+u)^2 - x^2 + 1}{2(1+u)} \Big) \right),
    \end{equation}
    and $L_k$ is the \emph{Legendre polynomial of degree} $k$.
\end{theorem}

\begin{proof}
    By Proposition \ref{prop:PDE}, the law $p_Z$ satisfies the weakly dissipative wave equation \eqref{dissipative-wave}. Setting $\tau \coloneqq 1 + u$ and using $\frac{a}{\theta}=k$, the PDE \eqref{dissipative-wave} can be written as the (shifted)
    EPD equation \eqref{EPD},
    which now starts at $\tau_0 \coloneqq 1$. 
    The Riemann-Volterra method, see \cite[Chapter 5]{Copson_1975} and Appendix \ref{appendix-proofs-A1}, applied to the EPD equation \eqref{EPD} yields the solution formula
    \begin{align}\label{EPD-solution}
        p(\tau,x) &= \frac{1}{2} \left(\frac{\tau_0}{\tau} \right)^k \big( p(\tau_0, x - \tau + \tau_0) +  p(\tau_0, x + \tau - \tau_0) \big) \\
        & + \frac{1}{2} \left(\frac{\tau_0}{\tau} \right)^k \int_{x - \tau + \tau_0}^{x + \tau - \tau_0} R(\tau, x; \tau_0, \xi) \partial_\tau p(\tau_0, \xi) \d \xi\\
        & + \frac{1}{2} \left(\frac{\tau_0}{\tau} \right)^k \int_{x - \tau + \tau_0}^{x + \tau - \tau_0} \Big( \frac{k}{\tau_0} R(\tau, x; \tau_0, \xi) - \partial_{\tau_0} R (\tau, x; \tau_0, \xi) \Big) p(\tau_0, \xi) \d \xi.
    \end{align}
    Here, $R(\tau, x; \tau_0, \xi)$ is the Riemann function corresponding to the EPD equation \eqref{EPD}, which can be given in terms of Gauss's hypergeometric function, see \eqref{Gauss-hypergeometric} in Appendix \ref{appendix-proofs-A1}. But since $k \in \N$, this simplifies to
    \begin{equation}\label{Riemann-function}
       R(\tau, x; \tau_0, \xi)
       =  L_{k-1}\left(\frac{\tau^2 + \tau_0^2 - (x-\xi)^2}{2\tau \tau_0} \right).
    \end{equation}
    Inserting $\tau_0 = 1$ and using the initial conditions $p(1, x) = \delta(x), \, \partial_{\tau}p(1, x) = 0$, the expression for $p(\tau,x)$ becomes
    \begin{align}\label{zwischenschritt}
        p(\tau,x) &= \frac{1}{2 \tau^k} \Big( \delta(x - \tau + 1)  + \delta(x + \tau - 1) \Big) \\
         &+ \frac{1}{2 \tau^k} \Big(k R(\tau,x;1, 0)  - \partial_{\tau_0} R (\tau,x;1, 0) \Big) \mathds{1}_{\{|x| < \tau-1\}}.
    \end{align}
    We are left to treat the absolutely continuous part. 
    First, we have
    \begin{equation}
        R(\tau,x;1, 0) =  L_{k-1}\left( \frac{\tau^2 - x^2 + 1}{2\tau} \right).
    \end{equation}
    Next, by the chain rule we compute
    \begin{align}
        \partial_{\tau_0} R (\tau,x;\tau_0, \xi) &= L_{k-1}'\left( \frac{\tau^2 + \tau_0^2 - (x-\xi)^2}{2\tau \tau_0} \right) \frac{4 (\tau_0)^2 \tau - 2\tau(\tau^2 + \tau_0^2 - (x-\xi)^2) }{(2\tau \tau_0)^2}.
    \end{align}
    Evaluated at $\tau_0 = 1$ and $\xi = 0$, this yields
    \begin{align}
        \partial_{\tau_0} R (\tau,x;1, 0) &= 
        L_{k-1}'\left( \frac{\tau^2 + 1 - x^2}{2\tau} \right) \frac{4 \tau - 2\tau(\tau^2 + 1 - x^2) }{(2\tau)^2}
        = L_{k-1}'\left( z \right) \Big(\frac{1}{\tau} - z \Big), 
    \end{align}
    where $z \coloneqq \frac{\tau^2 + 1 - x^2}{2\tau}$. Therefore, we obtain
    \begin{align}
        k R(\tau,x;1, 0) - \partial_{\tau_0} R (\tau,x;1, 0) 
        &= k  L_{k-1}\left( z \right) - L_{k-1}'\left( z \right) \Big(\frac{1}{\tau} - z \Big)\\
        &= k L_{k-1}(z) + z L_{k-1}'(z) - \frac{1}{\tau} L_{k-1}'(z).
    \end{align}
    Using the recurrence relation 
    \begin{equation}
        k L_{k-1}(z) + z L_{k-1}'(z) = L_k'(z),
    \end{equation}
    see \cite[Formula 28.21] {Spiegel2009},
    the absolutely continuous part in \eqref{zwischenschritt} becomes
    \begin{equation}
        \frac{1}{2 \tau^k} \Big( L_k'\left(\frac{\tau^2 + 1 - x^2}{2\tau} \right) - \frac{1}{\tau} L_{k-1}' \left(\frac{\tau^2 + 1 - x^2}{2\tau} \right) \Big) \mathds{1}_{\{|x| < \tau-1\}}.
    \end{equation}
    Finally, inserting $\tau = 1+u$ and noticing that $p_Z(u,x) = p(1+u, x)$ finishes the proof.
\end{proof}

\begin{remark}
    To study the law of the inhomogeneous Kac process $(Z_u)_u$, the weakly dissipative wave equation \eqref{dissipative-wave} was transformed to the EPD equation \eqref{EPD}, which -- \emph{importantly} -- starts at $\tau_0 = 1$. 
    Usual treatments of the EPD equation \cite{bresters1973, DeGregorioGarra2026, DeGregorioOrsingher2020, GarraOrsingher2016} consider the initial time $\tau_0 = 0$. Notice that there, the damping factor $\frac{2k}{\tau}$ admits a singularity at $\tau_0 = 0$. This changes significantly the qualitative behavior of the flow: for instance, the two Dirac singularities at the boundary of the wave do not emerge in that case. 
    However, in our case with $\tau_0 =1$, the singularity of the damping factor is avoided, and we obtain a rather simple and elegant representation of the law $p_Z$ via two Dirac points \eqref{legendre-formula-total} and two Legendre polynomials \eqref{legendre-formula-abs} for the cases $a = k \theta$.    
\end{remark}

\begin{remark}\label{remark:sophisticated}
    We want to note that for the general case $k = \frac{a}{\theta} \in \R$, the solution of the EPD equation \eqref{EPD} can still be expressed via the formula \eqref{zwischenschritt} involving the hypergeometric function \eqref{Gauss-hypergeometric}. 
    Yet, the simplification to \eqref{legendre-formula-abs} via Legendre polynomials crucially relies on the assumption $k = \frac{a}{\theta} \in \N$. 
\end{remark}

As a by-product, we derive the law of the OU-K process, the proof is given in Appendix \ref{appendix-proofs-A1}.

\begin{corollary}[Law of OU-K for $\frac{a}{\theta} \in \N$]\label{cor:law_OU-K}
The OU-K process $(Y_t)_t$ with speed $c > 0$, mean $\mu = 0$, switching intensity $a > 0$ and rate $\theta > 0$ fulfilling $\frac{a}{\theta} = k \in \N$ has the law 
    \begin{align}
        p_Y(t,x) 
        &= \frac{e^{\theta t}}{2(1+u)^k} \Big\{ p_{0} \left(e^{\theta t} (x - \frac{c(1-e^{-\theta t})}{\theta} ) \right) + p_{0} \left(e^{\theta t} (x + \frac{c(1-e^{-\theta t})}{\theta}) \right) \Big\} \\
        &+ \frac{\theta e^{2 \theta t}}{c} \int_{x - \frac{c(1-e^{-\theta t})}{\theta}}^{x + \frac{c(1-e^{-\theta t})}{\theta}} p_{0}(e^{\theta t}y) \, P_{Z}^{(k)}(u, \frac{\theta e^{\theta t}}{c}(x-y)) \d y,
    \end{align}
    where $p_0$ is the density of the law of $Y_0$, and $P_Z^{(k)}$ is given by the Legendre formula \eqref{legendre-formula-abs}.
\end{corollary}  

The representation of the OU-K law in the general case $\frac{a}{\theta} \in \R$ appears rather sophisticated involving hypergeometric functions, see Remark \ref{remark:sophisticated}, and also \cite[Section 5]{ratanov2021}.

\section{Inhomogeneous Kac and OU--K velocity field}\label{sec:velocity-kac}

Building on the explicit representation of the inhomogeneous Kac law in Theorem \ref{thm:main-polynomial-law}, we derive the (conditional) velocity fields associated with $(Z_u)_u$ and $(Y_t)_t$.
These velocity fields, which we prove to be \emph{uniformly bounded}, serve as the main ingredients of the OU-K generative flow matching model introduced in Section \ref{sec:generative-modeling}.
To the best of our knowledge, the velocity representations obtained in Theorem \ref{cor:velo-formula} and 
Proposition \ref{MSP-velo} are new to the literature.

\subsection{Velocity field of the inhomogeneous Kac process}\label{sbsec:vkac}

Before we explicitly compute the velocity field corresponding to $(Z_u)_u$, let us state an important (Lipschitz) regularity of our Kac processes. The proof follows the same lines as that in \cite[Proposition 5.6]{duong2026telegraphersgenerativemodelkac}, with only minor adjustments to the setting.

\begin{proposition}[Properties of the inhomogeneous Kac flow]\label{prop:lipschitz}
   Let $(Z_t)_t$ be a (possibly inhomogeneous) Kac process with arbitrary rate $a(t) > 0$ and speed $c>0$ defined in \eqref{inhom-Kac-process-Zt}.
   Then, for its probability distribution flow  $(p_t)_{t\ge 0} \coloneqq (P_{Z_t})_{t\ge 0}$, the following holds true:
   \begin{itemize}
    \item[\rm i)]
        $(p_t)_{t\ge 0}$ is Lipschitz continuous in the Wasserstein space $(\mathcal{P}_2(\R),W_2)$ with constant $c$, and $(Z_t)_t$ is almost surely Lipschitz continuous in $\R$ with the same Lipschitz constant.
        
     \item[\rm ii)] $(p_t)_{t\ge 0}$ is an absolutely continuous curve in $AC^2([0,T];\mathcal P_2(\R))$ for any $T>0$, and admits a velocity field $v(t,x)$ solving the continuity equation \eqref{eq:ce} and having the boundedness
     \begin{equation}\label{kac-boundedness}
    \|v_t\|_{L_2(p_t)}^2  \le  c^2 \quad \text{ for a.e. } t \in [0, \infty).
    \end{equation}
    
    \item[\rm iii)]
    The support of $(p_t)_{t\ge 0}$ grows at most linearly in time.
\end{itemize}
\end{proposition}

As a consequence of Theorem \ref{thm:main-polynomial-law}, we can explicitly calculate the corresponding velocity field of the probability curve induced by $(Z_u)_u$, given that $a$ is a multiple of $\theta$. The velocity formula \eqref{Z-velocity} will be essential to our flow matching approach in Section \ref{sec:generative-modeling}.

\begin{theorem}[Velocity field of the  inhomogeneous Kac process]\label{cor:velo-formula}
    Let $\frac{a}{\theta}=k \in \N$, and $(Z_u)_u$ be a symmetric, inhomogeneous Kac process with speed $c_Z = 1$ and switching rate $a(u) = \frac{a}{\theta(1+u)} = \frac{k}{1+u}$. 
    With $z\coloneqq z(u,x) \coloneqq \frac{(1+u)^2-x^2+1}{2(1+u)}$, define the velocity field
    \begin{align}\label{Z-velocity}
        v_Z(u,x) \coloneqq 
        \begin{cases}
            &\frac{x}{u + 1 - \frac{L_{k-1}'\left( z \right)}{L_{k}'\left( z \right)}}, \quad x \in (-u,u),\\
            &1, \hspace{20.7mm} x=u,\\
            &-1, \hspace{18.1mm} x=-u,\\
            &\text{arbitrary}, \hspace{9mm} \text{else}.
        \end{cases}
    \end{align}
       Then, the pair $(p_{Z}, v_Z)$ solves the continuity equation
    \begin{equation}\label{whole-CE}
        \partial_u \, p_{Z} + \partial_x (p_{Z} v_Z) = 0
    \end{equation}
    in the sense of distributions on $(u,x) \in (0,\infty)\times \R$.
\end{theorem}

\begin{proof}
    Define 
    \begin{equation}\label{flux-AC}
        J_Z^{\rm ac}(u,x) \coloneqq \frac{x}{2(1+u)^{k+1}} L_k' \left(z \right).
    \end{equation}
   1. First, we prove that $(p_Z, J_Z^{\rm ac})$ solves the continuity equation 
    \begin{equation}\label{interior-CE}
        \partial_u \, p_{Z} + \partial_x J_Z^{\rm ac} = 0 \quad 
        \text{for } x \in (-u,u).
    \end{equation}
    Let $\tau = 1 + u$. Then, $z = \frac{\tau^{2} + 1 - x^{2}}{2 \tau}$ and further
     \begin{align}\label{eq:derivativesflux}
    \partial_{u} z &= \frac{4 \tau^{2} - 2 (\tau^{2} + 1 - x^{2})}{4 \tau^{2}} = 1 - \frac{z}{\tau}  \qquad \text{and} \qquad \partial_{x} z = - \frac{x}{\tau}.
    \end{align}
    Let $x \in (-u,u)$. Differentiating \eqref{legendre-formula-abs} via the chain rule then yields
    \begin{align*}
    \partial_{u} P_Z^{(k)}(u, x) 
    &= \frac{1}{2 \tau^{k + 1}}\left(-k L_{k}'(z) + \frac{k + 1}{\tau}L_{k - 1}'(z) + (\tau - z)L_{k}''(z) - \left(1 - \frac{z}{\tau}\right)L_{k - 1}''(z)\right).
    \end{align*}
    Using the identities (\cite[Formulas 28.21, 28.22]{Spiegel2009}) 
    \begin{equation*}
    L_{k}'(\zeta) = k L_{k - 1}(\zeta) + \zeta L_{k - 1}'(\zeta) \qquad \text{and} \qquad L_{k - 1}'(\zeta) = \zeta L_{k}'(\zeta) - k L_{k}(\zeta),
    \end{equation*}
    we compute the derivatives
    \begin{align}\label{eq:firstthreeterm}
    L_{k}''(\zeta) \,&=\, k L_{k - 1}'(\zeta) + L_{k - 1}'(\zeta) + \zeta L_{k - 1}''(\zeta) 
    \,=\, (k + 1) L_{k - 1}'(\zeta) + \zeta L_{k - 1}''(\zeta),
    \\
    L_{k - 1}''(\zeta) \,&=\, L_{k}'(\zeta) + \zeta L_{k}''(\zeta) - k L_{k}'(\zeta) 
    \,=\, \zeta L_{k}''(\zeta) - (k - 1) L_{k}'(\zeta).
    \end{align}
    Plugging this into above expression for $\partial_{u} P_Z^{(k)}(u, x)$ gives
    \begin{align*}
    \partial_{u} P_Z^{(k)}(u, x) &= 
    \frac{1}{2 \tau^{k + 1} } \left( 
    -k L_{k}'(z) + \frac{1}{\tau}(L_{k}''(z) - z L_{k - 1}''(z)) + (\tau - z)L_{k}''(z) - \left(1 - \frac{z}{\tau}\right)L_{k - 1}''(z) \right) \\
    &= \frac{1}{2 \tau^{k + 1} } \left( -k L_{k}'(z) + \left(\frac{1}{\tau} + \tau - z\right)L_{k}''(z) - L_{k - 1}''(z) \right) \\ 
    &= \frac{1}{2 \tau^{k + 1} } \left( -k L_{k}'(z) + \left(\frac{1}{\tau} + \tau - z\right)L_{k}''(z) - z L_{k}''(z) + (k - 1)L_{k}'(z) \right)\\
    &= \frac{1}{2 \tau^{k + 1} } \left( -L_{k}'(z) + \left(\frac{1}{\tau} + \tau - 2 z\right)L_{k}''(z) \right).
    \end{align*}
    Since
    \begin{equation*}
        2 z = \frac{\tau^{2} + 1 - x^{2}}{\tau} = \tau + \frac{1}{\tau} - \frac{x^{2}}{\tau},
    \end{equation*}
    we obtain
    \begin{equation*}
        \partial_{u} P_Z^{(k)}(u, x) = \frac{1}{2 \tau^{k + 1}} \left(-L_{k}'(z) + \frac{x^{2}}{\tau}L_{k}''(z)\right).
    \end{equation*}
    At the same time, differentiating \eqref{flux-AC} with respect to $x$ results via \eqref{eq:derivativesflux} in
    \begin{align*}
    \partial_{x} J_Z^{\rm ac}(u,x) &= \frac{1}{2 \tau^{k + 1}} \left( L_{k}'(z) - \frac{x^{2}}{\tau}L_{k}''(z)\right).
    \end{align*}
    Therefore, it holds the continuity equation \eqref{interior-CE} for interior points $x \in (-u,u)$.
\\[1ex]
2.    Arguing in the same way as in the proof of \cite[Theorem 4.4]{duong2026telegraphersgenerativemodelkac}, Proposition \ref{prop:lipschitz} gives the \emph{existence} of some velocity field $v^*(u,x)$, which is antisymmetric in $x$ and solves the continuity equation \eqref{whole-CE} on the whole $(0,\infty)\times \R$. We need to show that $v^*(u,x)$ coincides with \eqref{Z-velocity}. For a.e.\ $u > 0$,
    the symmetry of $P_Z^{(k)}(u, \cdot)$ in \eqref{legendre-formula-abs}, the antisymmetry of $J_Z^{\rm ac}(u,\cdot)$ in \eqref{flux-AC}, and the first step \eqref{interior-CE} yield 
    \begin{equation}
       2 \, P_Z^{(k)}(u, x) \, v^*(u,x) = \int_{-x}^x \partial_x \big(P_Z^{(k)} v^* \big) \d x
       = - \int_{-x}^x \partial_u P_Z^{(k)} \d x
       = \int_{-x}^x \partial_x J_Z^{\rm ac} \d x
       = 2 J_Z^{\rm ac}(u,x),
    \end{equation}
    for almost all $x \in (-u,u)$, and hence,
    \begin{align}
        v^*(u, x) &= \frac{J_Z^{\rm ac}(u, x)}{P_Z^{(k)}(u, x)} = \frac{\frac{x}{2(1 + u)^{k + 1}}L_{k}'(z)}{\frac{1}{2(1 + u)^{k}}\left(L_{k}'(z) - \frac{1}{1 + u}L_{k - 1}'(z)\right)} \\
        &= \frac{x L_{k}'(z)}{(1 + u) L_{k}'(z) - L_{k - 1}'(z)} \\ 
        &= \frac{x}{1 + u - \frac{L_{k - 1}'(z)}{L_{k}'(z)}}
    \end{align}
    for almost all $x \in (-u,u)$, which proves that $v^*_u$ coincides with \eqref{Z-velocity} for a.e.\ interior point $x \in (-u,u)$.
    \\[1ex]
    3. We need to determine the velocity field $v^*(u,x)$ on the boundary points $x = \pm c u$, where we keep track on the speed constant $c \coloneqq c_Z \coloneqq 1$.
For this, let $\phi \in C_c^\infty(\R)$ be any test function. Then, we get
\begin{align}
    \int_{-\infty}^\infty \partial_u p_Z(u,x) \phi(x) \d x
    &= - \frac{k}{2(1+u)^{k+1}} \big(\phi(c u) + \phi(-c u) \big) + \frac{1}{2(1+u)^k} \big(-c \phi'(-c u) + c \phi'(c u) \big) \\
    & \quad +\int_{-c u}^{c u} \partial_u P_Z^{(k)} (u,x) \phi(x) \d x.
\end{align}
On the other hand, since $(p_Z, v^*)$ solves the continuity equation \eqref{whole-CE}, it holds for a.e.\ $u>0$ that
\begin{align}
    &\int_{-\infty}^\infty \partial_u p_Z(u,x) \phi(x) \d x
    = \int_{-\infty}^\infty - \partial_x (p_Z(u,x) v^*(u,x)) \phi(x) \d x
    = \int_{-\infty}^\infty p_Z(u,x) v^*(u,x) \phi'(x) \d x \\
    &= \frac{1}{2(1+u)^k} \big(v^*(u, -c u) \phi'(-c u) + v^*(u, c u) \phi'(c u) \big) + 
    \int_{-c u}^{c u} P_Z^{(k)}(u,x) v^*(u,x) \phi'(x) \d x.
\end{align}
Now, choose a sequence of test functions $\phi_n \in C_c^\infty(\R)$ such that $\phi_n \equiv 0$ on $[-c u + \frac{1}{n}, c u - \frac{1}{n}]$, $\|\phi_n\|_\infty < \frac{1}{n}$, $\phi_n'(c u) = 1$ and $\phi_n'(-c u) = -1$. Then, letting $n \to \infty$ in the above identities, it follows
\begin{align}
  \frac{1}{2(1+u)^k}  \big(-c \cdot (-1) + c \cdot 1 \big) = \frac{1}{2(1+u)^k} \big(v^*(u, -c u) \cdot (-1) + v^*(u,c u) \cdot 1 \big).
\end{align}
Finally, since $c = 1$, the antisymmetry of $v^*(u,\cdot)$ delivers $v^*(u,u) = 1$ and $v^*(u,-u) = -1$, which shows that $v^*$ has the form \eqref{Z-velocity} for $p_u^Z$-a.e.\ $x \in \R$.
\end{proof}

\begin{remark}
    It can be shown that $1 - \frac{L_{k-1}'\left( z \right)}{L_{k}'\left( z \right)} > 0$ for all $u>0$ and $x \in (-u,u)$, see Lemma \ref{lem:legendre-quotient-estimate} below. Therefore, we even have the pointwise bound
    \begin{equation}\label{kac-boundedness-pointwise}
        |v_Z(u,x)| \le 1 = c_Z  \quad \text{for all } u \ge 0, ~ x \in \R,
    \end{equation}
    that is, the velocity field $v_Z$ is globally bounded by the Kac process' speed $c_Z$, and not only in the $\|\cdot\|_{L_2(p_t)}$-norm as in \eqref{kac-boundedness}.
\end{remark}

\subsection{Conditional velocity field of the OU--K process}\label{sbsec:vkac-cond}

In this section, we derive the conditional velocity fields appearing in the CFM loss \eqref{cfm_loss}, forming the basis of our OU-K generative model proposed in Section \ref{sec:generative-modeling}.

Let us start with the characterisation of the velocity fields of Markov switching processes. 
\begin{proposition}[Velocity field of Markov switching processes]\label{MSP-velo}
    Let $(X_t, I_t)_t$ be a Markov switching process with characteristic tuple $(F, \lambda)$ and $\mathbb{E}[|X_0|^2] < \infty$. Let $p_t$ be the law of $X_t$ for all $t \geq 0$. Then the curve $(p_t)_t$ is an absolutely continuous curve in the Wasserstein space $\mathcal{P}_2(\R^d)$, and 
    \begin{equation}\label{eq:MSP-velo}
        v_t(x) \coloneqq \mathbb{E}[ F^{I_t}(X_t) \mid X_t = x]
    \end{equation}
    defines a velocity field associated with $(p_t)_t$.
\end{proposition}
    The result extends the well-known time-derivative formula for velocity fields, see \cite[Theorem 6.3]{WS2025}, since by definition \eqref{eq:mspdef} it holds 
    $$\mathbb{E}[F^{I_t}(X_t) \mid X_t = x] = \mathbb{E}[\dot X_t \mid X_t = x].$$ 
    Note that this only requires piecewise regularity of the Markov switching processes.
  
\begin{proof}[{Proof of Proposition \ref{MSP-velo}}]

We define the vector-valued flux 
\begin{equation}
    J_t \coloneq \sum_{i = 1}^n F^i \rho_{i, t}.
\end{equation}
Summing up the mode-wise Kolmogorov forward equations \eqref{mode-wise-KFE} over $i=1,\dots,n$, and using that the mode-wise laws $\rho_{i, t}$ sum up to $p_t$, we see that the switching terms cancel and we obtain
\begin{equation}\label{CE-in-proof}
    \partial_t p_t = \sum_{i = 1}^n \left(- \operatorname{div}(F^i\rho_{i, t}) + \sum_{j \neq i} (\lambda^{ji}\rho_{j, t} - \lambda^{ij}\rho_{i, t})\right)  = - \operatorname{div} \left( \sum_{i = 1}^n F^i \rho_{i, t} \right) = -\operatorname{div}(J_t)
\end{equation}
in the sense of distributions on $(0,\infty) \times \R^d$. As $\rho_{i, t} \ll p_t$, it follows $J_t \ll p_t$. Computing the Radon-Nikodym derivative
\begin{equation}\label{msp-velocity-formula}
    v_t(x) \coloneq \frac{\d J_t}{\d p_t}(x) = \sum_{i = 1}^n F^i(x) \frac{\d \rho_{i, t}}{\d p_t}(x) = \mathbb{E}[F^{I_t}(X_t) \mid X_t = x],
\end{equation}
it holds $v_t p_t = J_t$, and \eqref{CE-in-proof} yields
\begin{equation}\label{eq:MSP-velo-CE}
    \partial_t p_t + \operatorname{div}(v_t p_t) = 0
\end{equation}
in the sense of distributions on $(0,\infty) \times \R^d$. \par 

The absolute continuity of $(p_t)_t$ in $\mathcal{P}_2(\R^d)$ follows by standard arguments, we defer the details to Appendix \ref{appendix-proofs-A2}.
\end{proof}

Proposition \ref{MSP-velo} allows us to derive the conditional velocity field of an OU-K process $(Y_t)_t$ by connecting it to the velocity of the inhomogeneous Kac process $(Z_u)_u$ of Proposition \ref{prop:transformation}.

\begin{proposition}[Conditional OU-K velocity]\label{pro:OU-K-velo}
Assume that the counting process $(N_t)_t$ (and thus, also $(I_t)_t$, $(Z_u)_u$) is independent of the data variable $Y_0$,
and let $(Y_t)_t$ be the OU-K process with speed $c > 0$, mean $\mu = 0$, switching intensity $a > 0$ and rate $\theta > 0$. \\
Then, the process $(Y_t)_t$ \emph{conditioned} to start in $Y_0 = y_0$
has the associated  conditional velocity field
\begin{align}\label{Y-cond-velo}
    v_t^Y(y|y_0) 
    &= -\theta y + c \, v_Z \big(e^{\theta t}-1, ~ \frac{\theta}{c}( e^{\theta t}y - y_0) \big),
\end{align}
with a velocity field $v_Z$ associated with $(Z_u)_u$.
\end{proposition}

\begin{proof}
First, since $(Y_t, I_t)$ defines a Markov switching process with tuple $(F,\lambda)$ given in Proposition \ref{Z-Kac-is-MSP}, it holds by Proposition \ref{MSP-velo} that
\begin{align}
    v_t^Y(y|y_0) 
    &= \mathbb{E}[F^{I_t}(Y_t) \mid  Y_t = y, Y_0 = y_0]\\
    &= -\theta y + c \, \mathbb{E}[ (-1)^{I_t} \mid Y_t = y, ~ Y_0 = y_0].
\end{align}
With the corresponding inhomogeneous Kac process $(Z_u)_u$ with speed $c_Z = 1$ and rate $a(u) = \frac{a}{\theta(1+u)}$ from Proposition \ref{prop:transformation}, it holds via \eqref{relation-to-Z_t} that
\begin{align}
    v_t^Y(y|y_0) &= -\theta y + c \, \mathbb{E}[ (-1)^{I_t} \mid e^{-\theta t}x_0 + \frac{e^{-\theta t}}{\theta} c Z_{e^{\theta t}-1} = y, ~ Y_0 = y_0]\\
    &= -\theta y + c \, \mathbb{E}[ (-1)^{I_t} \mid Z_{e^{\theta t}-1} = \frac{\theta}{c}( e^{\theta t}y - y_0)],
\end{align}
where we took the independence of $I_t$ and $Y_0$ into account. Next, notice that the Poisson process $(N_t)_t$ defining $(Y_t)_t$ via \eqref{def:OU-K} can be written in terms of the inhomogeneous Poisson process $(\tilde N_u)_u$ from \eqref{inhom-Kac-process-Zt} with rate $a(u) = \frac{a}{\theta(1+u)}$ via
\begin{equation}
    N_t = \tilde N_{e^{\theta t}-1}.
\end{equation}
This yields 
\begin{equation}
    (-1)^{I_t} = B_\frac{1}{2} (-1)^{N_t} = B_\frac{1}{2} (-1)^{\tilde N_{e^{\theta t}-1}}
    = \partial_u Z_u |_{u =e^{\theta t}-1}.
\end{equation}
Finally, we arrive at
\begin{align}
    v_t^Y(y|y_0) 
    &= -\theta y + c \, \mathbb{E}[ \partial_u Z_u |_{u =e^{\theta t}-1} \mid Z_{e^{\theta t}-1} = \frac{\theta}{c}( e^{\theta t}y - y_0)]\\
    &= -\theta y + c \, v_Z \big(e^{\theta t}-1, ~ \frac{\theta}{c}( e^{\theta t}y - y_0) \big).
\end{align}
\end{proof}

Hence, given $y \sim Y_t$ and $y_0 \sim Y_0$, we only need to calculate the velocity field $v_Z(u, x)$ at the transformed time $u = e^{\theta t}-1$ and position $x = \frac{\theta}{c}( e^{\theta t}y - y_0)$. In the case $a = k\theta$, the velocity field $v_Z(u, x)$ is simply given by the Legendre formula \eqref{Z-velocity}, also see Appendix \ref{app:odd-flux}.

\subsection{Extension to \texorpdfstring{$\R^d$}{Rd} and time schedules}\label{subsec:multi-d}
Since our above results were considered in the one-dimensional case, we now describe the transition to a setting suitable for generative modeling.

\paragraph{Multi-Dimensions.}
Note that by the conditioning on fixed points $Y_0 = y_0$, 
the above procedure \eqref{Y-cond-velo} is still valid if we consider a vector-valued process $\mathbf{Z}_u = (Z_u^1, \dots, Z_u^d)$ consisting of independent, inhomogeneous Kac processes $(Z_u^i)_{u}$ as follows.

\begin{proposition}[Conditional OU-K velocity in $\R^d$]\label{pro:OU-K-velo-Rd}
Let $\mathbf{Y}_0$ be a random variable in $\R^d$, and let $(\mathbf{Y}_t)_t$ be the multi-dimensional OU-K process given as in \eqref{relation-to-Z_t},
with speed $c > 0$, mean $\mu = 0$, switching intensity $a > 0$ and rate $\theta > 0$;
and a corresponding vector-valued process $\mathbf{Z}_u = (Z_u^1, \dots, Z_u^d)$ consisting of independent, inhomogeneous Kac processes $(Z_u^i)_{u \ge 0}$.

Assume that the counting process $(\mathbf{N}_t)_t$ (and thus, also $(\mathbf{I}_t)_t$, $(\mathbf{Z}_u)_u$) is independent of the data variable $\mathbf{Y}_0$.
Then, the process $(\mathbf{Y}_t)_t$ \emph{conditioned} to start in $\mathbf{Y}_0 = y_0 \in \R^d$
has the associated  conditional velocity field
\begin{align}\label{Y-cond-velo-Rd}
    v_t^{\mathbf{Y}}(y|y_0) 
    &= -\theta y + c \, v_{\mathbf{Z}} \big(e^{\theta t}-1, ~ \frac{\theta}{c}( e^{\theta t}y - y_0) \big),
\end{align}
with a velocity field $v_{\mathbf{Z}}$ admitting the decomposition
\begin{equation}\label{cond-velo-multi-d}
    v_{\mathbf{Z}}(u, x_1,\dots, x_d) = (v_{Z_1}(u, x_1), \dots, v_{Z_d}(u, x_d)),
\end{equation}
where the one-dimensional vector fields $v_{Z_i}$ are associated with $(Z_u^i)_u$, $i=1,\dots,d$.
\end{proposition}
\begin{proof}
    Following the exact lines of the proof of Proposition \ref{pro:OU-K-velo}, we arrive at
\begin{align}
    v_t^{\mathbf{Y}}(y|y_0) 
    &= -\theta y + c \, \mathbb{E}[ \partial_u \mathbf{Z}_u |_{u =e^{\theta t}-1} \mid \mathbf{Z}_{e^{\theta t}-1} = \frac{\theta}{c}( e^{\theta t}y - y_0)].
\end{align}
Now, setting $\tilde u \coloneqq e^{\theta t}-1$ and $x \coloneqq \frac{\theta}{c}( e^{\theta t}y - y_0) \in \R^d$, we infer by the independence of the components $Z_u^i$, $i=1,\dots,d$, that
\begin{align}
    \mathbb{E}[ \partial_u \mathbf{Z}_u |_{u =\tilde u} \mid \mathbf{Z}_{\tilde u} = x]
    &= \left(  \mathbb{E}[ \partial_u Z_u^1 |_{u =\tilde u} \mid \mathbf{Z}_{\tilde u} = x], \, \dots, \, 
    \mathbb{E}[ \partial_u Z_u^d |_{u =\tilde u} \mid \mathbf{Z}_{\tilde u} = x]\right)\\
    &= \left(  \mathbb{E}[ \partial_u Z_u^1 |_{u =\tilde u} \mid Z^1_{\tilde u} = x_1], \, \dots, \,
    \mathbb{E}[ \partial_u Z_u^d |_{u =\tilde u} \mid Z^d_{\tilde u} = x_d]\right)\\
    &\eqcolon (v_{Z_1}(\tilde u, x_1), \dots, v_{Z_d}(\tilde u, x_d)),
\end{align}
where $v_{Z_i}$ is associated with $(Z_u^i)_u$ by Proposition \ref{MSP-velo}.
\end{proof}

The decomposition \eqref{cond-velo-multi-d} holds more generally for velocity fields of factorizing probability curves, see also \cite[Lemma 6.1]{duong2026telegraphersgenerativemodelkac}. 
We emphasize that the componentwise independence assumption applies \emph{exclusively}
to the noising process $(\mathbf{Z}_u)_u$, and imposes \textbf{no restrictions whatsoever on the data variable} $\mathbf{Y}_0$. Withing flow matching, the cross-dimensional dependencies of $\mathbf{Y}_0$ are reconstructed through the averaging of conditional velocity fields in the CFM loss \eqref{cfm_loss}.

This allows us to apply the above one-dimensional constructions in Section \ref{sec:generative-modeling}, where we consider high-dimensional generative modeling problems.

\paragraph{Time schedules.}

Motivated by the following observation, we want to consider a rescheduled version of the OU-K process.
\begin{remark}\label{remark:time-schedule}
   Recall that choosing $\theta = \frac{1}{2}$, $\mu=0$, and $\sigma = 1$ in \eqref{def:OU} yields the prototypical OU process
   \begin{equation}
       X_t \overset{d}{=} 
       e^{-\frac{t}{2}} X_0 + \sqrt{1 - e^{-t}} \, \Xi.
   \end{equation}
   This process is applied in score-based diffusion models under various time schedules given in the form $g(t) = \int_0^t \beta(s) \d s, ~ t \in [0,1)$, see \cite{Song2021, nichol2021improved}. That is, the process $\tilde X_t \coloneqq X_{g(t)}$, $t \in [0,1)$, is considered, satisfying the SDE
   \begin{equation}\label{scheduled-OU-dynamics}
       \d \tilde X_t = - \frac{1}{2} \beta(t) \tilde X_t \d t + \sqrt{\beta(t)} \d W_t, \quad t \in (0,1).
   \end{equation}
   \end{remark}
This motivates us to apply our Ornstein-Uhlenbeck-Kac process within a generative model in the \emph{rescheduled} form 
   \begin{equation}
       \tilde Y_t \coloneqq Y_{g(t)}, \quad t \in [0,1),
   \end{equation}
where $(Y_t)_{t \in [0,\infty)}$ is the OU-K process from \eqref{def:OU-K} with $\mu=0$, $\theta \in \{\frac{1}{2}, 1\}$, suitably chosen $a,c > 0$, 
and $g:[0,1)\to [0,\infty)$ is some sufficiently smooth time schedule with $g(0)=0$ and $\lim_{t \uparrow 1} g(t) = \infty$.

The corresponding conditional velocity field of $(\tilde Y_t)_{t \in [0,1)}$ can be obtained via the chain rule as 
\begin{align}\label{cond-velo-schedule}
   v_t^{\tilde Y}(y|y_0) &= \mathbb{E}[ \partial_t (Y_{g(t)}) \mid Y_{g(t)} = y, \, Y_{g(0)} = y_0]\\
   &=  g'(t) \, \mathbb{E}[ \partial_t Y_{t}|_{t = g(t)} \mid Y_{g(t)} = y, \, Y_0 = y_0]\\
   &= g'(t) \, v_{g(t)}^{Y}(y | y_0),
\end{align}
where the latter velocity field $v^Y$ can be traced back to \eqref{Y-cond-velo}.

The conditional velocity field \eqref{cond-velo-schedule} can be shown to be \emph{bounded} under weak assumptions on the schedule $g$. This is the content of the following section.

\paragraph{Compressed OU-K process.} 
Recall the representation obtained in Proposition \ref{prop:transformation} given as
\begin{equation}\label{1d-backbone}
   \tilde Y_{t} = Y_{g(t)} = e^{- \theta g(t)} Y_{0} + \frac{c}{\theta}e^{-\theta g(t)} Z_{e^{\theta g(t)} - 1},
\end{equation}
where $(Z_u)_{u \in [0, \infty)}$ is a unit-speed, inhomogeneous Kac process with switching rate given in \eqref{specific-switching-OU-K}. Setting $\tau \coloneqq \tau(t) \coloneqq 1 - e^{-\theta g(t)}, \,  t \in [0,1)$, it holds $\tau(0)=0$, $\lim_{t\uparrow 1}\tau(t) = 1$, and 
\begin{equation}\label{physical-process-being-compressed}
    \tilde Y_t = (1 - \tau)Y_0 + \frac{c}{\theta} (1 - \tau) Z_{\frac{\tau}{1 - \tau}}, \quad t \in [0,1).
\end{equation}
Defining the \emph{(time-)compressed} versions of the inhomogeneous Kac and OU-K process as
\begin{equation}\label{compressed-processes}
    \hat Z_\tau \coloneq (1 - \tau)Z_{\frac{\tau}{1 - \tau}}, \quad 
    \hat Y_{\tau} \coloneq (1 - \tau)Y_0 + \frac{c}{\theta} \hat Z_\tau, \qquad \tau \in [0, 1),
\end{equation}
it holds 
\begin{equation}\label{rescheduled-relations}
    \tilde Y_{t} = Y_{g(t)} = \hat Y_{\tau(t)} \quad \text{for all } t \in [0, 1).
\end{equation}
Next, we study the compressed processes $\hat Z_{\tau \in [0,1)}$ and $\hat Y_{\tau \in [0,1)}$ and their associated velocity fields. To this end, we have the following boundedness property.

\begin{proposition} [Boundedness of the compressed Kac velocity field]\label{prop:compressed-Z-velo}
    Let $(Z_u)_{u \geq 0}$ be the inhomogeneous Kac process with associated velocity $v_Z$ from Theorem \ref{cor:velo-formula}, and $(\hat Z_\tau)_{\tau \in [0, 1)}$ be its compression. Define the velocity field
    \begin{align}\label{compressed-Z-velocity}
        v_{\hat Z}(\tau, x) \coloneqq \frac{1}{1 - \tau} \left(v_Z\left(\frac{\tau}{1 - \tau}, \frac{x}{1 - \tau}\right) - x \right), \quad \tau \in (0,1), \, x \in [-\tau,\tau].
    \end{align}
    Then, $v_{\hat Z}$ is a velocity field associated with $(\hat Z_\tau)_{\tau \in [0, 1)}$, and on its support $D \coloneq \{ (\tau, x) \mid \tau \in (0, 1), \, x \in [-\tau, \tau]\}$, it has the representation
    \begin{align}\label{compressed-velocity-legendre}
        v_{\hat Z}(\tau, x) = \frac{x\frac{L_{k - 1}'(z)}{L_k'(z)}}{1 - (1 - \tau)\frac{L_{k - 1}'(z)}{L_k'(z)}} \quad \text{ for }\, |x| < \tau,
        \qquad v_{\hat Z}(\tau, \pm \tau) = \pm 1,
    \end{align}      
    where $z \coloneqq \frac{1 - x^2 + (1 - \tau)^2}{2(1 - \tau)}$, and satisfies the uniform bound 
    \begin{equation}\label{compressed-Z-velo-bound}
        \sup_{(\tau, x) \in D} |v_{\hat Z}(\tau, x)| \leq 1.
    \end{equation}
\end{proposition}

Note that the singularity for $\tau \to 1$ in \eqref{compressed-Z-velocity} is removable by means of the representation \eqref{compressed-velocity-legendre}. Indeed, the velocity field $v_{\hat Z}$ remains bounded as $\tau \to 1$ and $\tau \to 0$. 

The proof of Proposition \ref{prop:compressed-Z-velo} makes use of the following lemma, which bounds the quotient of the Legendre polynomials appearing in \eqref{Z-velocity}. We provide its proof in Appendix \ref{appendix-proofs-A2}.

\begin{lemma}\label{lem:legendre-quotient-estimate}
    Let $k \in \N$, $k \ge 1$, and $L_k, L_{k-1}$ be the Legendre polynomials of degree $k$ and $k-1$, respectively. Then it holds $L_k'(z) > 0$ for all $z\ge 1$, and further,
    \begin{equation}\label{legendre-quotient}
        \left| \frac{L_{k - 1}'(z)}{L_k'(z)}\right| 
        = \frac{L_{k - 1}'(z)}{L_k'(z)} < 1 \quad \text{for all } z \ge 1.
    \end{equation}
\end{lemma}

\begin{proof}[{Proof of Proposition \ref{prop:compressed-Z-velo}}] 
    Since $\partial_\tau \frac{\tau}{1 - \tau} = \frac{1}{(1 - \tau)^2}$, we compute the time-derivative
    \begin{equation}\label{derivative-of-compressed-inhom-kac}
        \partial_\tau \hat Z_\tau = - Z_{\frac{\tau}{1 - \tau}} + \frac{1}{1 - \tau} \partial_u Z_{u}\mid_{u = \frac{\tau}{1 - \tau}} \qquad \text{for almost all } \tau \in [0, 1).
    \end{equation}
    Taking the conditional expectation of \eqref{derivative-of-compressed-inhom-kac} provides a velocity field associated with $(\hat Z_\tau)_\tau$, \emph{cf}.\ Proposition \ref{MSP-velo}. 
    As $\hat Z_\tau = x \Leftrightarrow Z_{\frac{\tau}{1 - \tau}} = \frac{x}{1 - \tau}$, this reads as
    \begin{align}\label{compressed-inhomo-kac-velocity}
        v_{\hat Z}(\tau, x) &= \mathbb{E}[\partial_\tau \hat Z_\tau \mid \hat Z_\tau = x]  \\ 
        &= - \frac{x}{1 - \tau} + \frac{1}{1 - \tau} \mathbb{E}\left[\partial_u Z_u \mid_{u = \frac{\tau}{1 - \tau}} \mid Z_{\frac{\tau}{1 - \tau}} = \frac{x}{1 - \tau}\right] \\ 
        &= \frac{1}{1 - \tau} \left(- x + v_Z\left(\frac{\tau}{1 - \tau}, \frac{x}{1 - \tau}\right) \right),
    \end{align}
    where $v_Z$ is the velocity field \eqref{Z-velocity}.   
    Now, it immediately follows $v_{\hat Z}(\tau, \pm \tau) = \pm 1$.
    Next, let $\tau \in (0,1)$ and $|x| < \tau$. The argument $z$ appearing in the Legendre polynomials is given as
    \begin{equation}
        z \coloneqq \frac{\left(1 + \frac{\tau}{1 - \tau}\right)^2 - \left(\frac{x}{1 - \tau}\right)^2 + 1}{2\left(1 + \frac{\tau}{1 - \tau}\right)}=  \frac{\frac{1}{(1 - \tau)^2}}{\frac{1}{1 - \tau}} \frac{1 - x^2 + (1 - \tau)^2}{2} = \frac{1 - x^2 + (1 - \tau)^2}{2(1 - \tau)} > 1.
    \end{equation}
    Using \eqref{Z-velocity}, we obtain
    \begin{equation}
        v_Z\left(\frac{\tau}{1 - \tau}, \frac{x}{1 - \tau}\right) = \frac{\frac{x}{1 - \tau}}{1 + \frac{\tau}{1 - \tau} - \frac{L_{k - 1}'(z)}{L_k'(z)}}  = \frac{x}{1 - (1 - \tau) \frac{L_{k - 1}'(z)}{L_{k}'(z)}},
    \end{equation}
    which proves the formula \eqref{compressed-velocity-legendre} via
    \begin{align}
        v_{\hat Z}(\tau, x) = \frac{1}{1 - \tau} \left(\frac{x}{1 - (1 - \tau) \frac{L_{k - 1}'(z)}{L_k'(z)}} - x \right) = \frac{x\frac{L_{k - 1}'(z)}{L_k'(z)}}{1 - (1 - \tau)\frac{L_{k - 1}'(z)}{L_k'(z)}}. 
    \end{align}
    Notice that the singularity for $\tau \to 1$ has been removed.
    By Lemma \ref{lem:legendre-quotient-estimate}, it holds that $\left| \frac{L_{k - 1}'(z)}{L_k'(z)}\right| < 1$ for all $z \geq 1$. Then, for $|x| < \tau < 1$, \eqref{compressed-velocity-legendre} admits the bound 
    \begin{equation}
        |v_{\hat Z}(\tau, x)| <  \frac{\tau}{1-(1-\tau) \left| \frac{L_{k - 1}'(z)}{L_k'(z)}\right|} < \frac{\tau}{1-(1-\tau)} = 1.
    \end{equation}
\end{proof}

With Proposition \ref{prop:compressed-Z-velo}, we obtain the boundedness also of the conditional velocity field of the compressed OU-K process.
\begin{corollary}\label{cor:conditional-bound}
    Let $(\hat Y_\tau)_{\tau \in [0, 1)}$ be a compressed OU-K process with rate $\theta > 0$ and speed $c > 0$. Then, its conditional velocity field is given by
    \begin{equation}\label{compressed-velocity}
        v_\tau^{\hat Y}(y \mid y_0) = - y_0 + \frac{c}{\theta}v_{\hat Z}(\tau, x),
    \end{equation}
    where $x \coloneqq  \frac{\theta}{c}(y - (1 - \tau)y_0)$, and $v_{\hat Z}$ is the velocity field \eqref{compressed-Z-velocity}. Further, the conditional velocity admits the uniform bound
    \begin{equation}\label{compressed-conditional-velo-bound}
        |v^{\hat Y}_\tau (y \mid y_0)|  \leq |y_0| + \frac{c}{\theta}  \quad \text{for all } \tau \in [0, 1).
    \end{equation}
\end{corollary}

\begin{proof}
    Differentiating the representation \eqref{compressed-processes} of $\hat Y_\tau$ and taking the conditional expectation yields
    \begin{equation}\label{conditional-compressed-OU-velo}
        v^{\hat Y}_\tau(y \mid y_0) 
        = -y_0 + \frac{c}{\theta} v_{\hat Z}(\tau, x).
    \end{equation}
    Now, the bound \eqref{compressed-conditional-velo-bound} follows immediately by \eqref{compressed-Z-velo-bound}.
\end{proof}

The above result states a stark contrast between the OU-K process and the usual OU process:
\begin{remark}\label{rem:difference-OU-K-OU}
Analogously to the compressed OU-K process, we can define a compressed version of the classical OU diffusion process satisfying the representation
\begin{equation}
    \hat X_\tau \overset{d}{=} (1 - \tau)X_0 + \frac{\sigma}{\sqrt{2 \theta}} \sqrt{(2 - \tau)\tau} \, \Xi, \qquad \tau \in [0, 1),
\end{equation}
where $\Xi \sim \mathcal{N}(0, I)$.
Using \cite[Theorem 3]{lipman2023flow}, (emphasizing our data $\to$ noise convention) we obtain the conditional velocity
\begin{equation}
    v_\tau^{\hat X}(x \mid x_0) \,=\, -x_0 + \frac{1 - \tau}{(2 - \tau)\tau}(x - (1 - \tau)x_0) \,=\, \frac{(1-\tau) x - x_0}{(2 - \tau)\tau}, \qquad \tau \in [0, 1).
\end{equation}
Notice that this velocity field admits the blow-ups $\lim_{|x|\to \infty} |v_\tau^{\hat X}(x \mid x_0)| = \infty$, and more importantly, $\lim_{\tau \downarrow 0} |v_\tau^{\hat X}(x \mid x_0)| = \infty$ in the limit $\tau \downarrow 0$ approaching the data, in contrast to the OU-K bound \eqref{compressed-conditional-velo-bound}.
The latter blow-up also occurs for the rescheduled OU process, see Appendix \ref{appendix-proofs-A2}.
\end{remark}

As a by-product of Proposition \ref{prop:compressed-Z-velo}, we obtain the Lipschitz regularity in Wasserstein-2 of the compressed Kac probability flows, which allows them to be naturally extended to the final time $\tau = 1$. The proof is deferred to Appendix \ref{appendix-proofs-A2}.

\begin{proposition}[Wasserstein extension of the compressed OU-K law]\label{prop:compressed-W2-extension}
    Let $(\hat Z_\tau)_{\tau \in [0, 1)}$ and $(\hat Y_\tau)_{\tau \in [0, 1)}$ be the compressed processes from \eqref{compressed-processes} with probability curves $(p^{\hat Z}_\tau)_{\tau \in [0, 1)}$ and $(p^{\hat Y}_\tau)_{\tau \in [0, 1)}$. Assume $\mathbb{E}[|Y_0|^2] < \infty$ and let $p_\infty$ be the invariant distribution \eqref{invariant-OU-K} of $(Y_t)_{t \geq 0}$. Then the following holds:
   \begin{enumerate}[label=(\roman*)]
       \item[\rm i)] Setting $p_1^{\hat Z} \coloneq \left(T_\frac{\theta}{c}\right)_{\#} p_\infty$, where $T_\frac{\theta}{c}:\R \to \R,\ x \mapsto \frac{\theta}{c}x$, the curve $(p^{\hat Z}_\tau)_{\tau \in [0, 1]}$ is Lipschitz continuous in Wasserstein-2, i.e.\ it holds       
       \begin{equation}
           W_2(p_{\tau_1}^{\hat Z}, p_{\tau_2}^{\hat Z}) \leq |\tau_1 - \tau_2| \quad \text{for all }\, \tau_1, \tau_2 \in [0, 1].
       \end{equation}
       \item[\rm ii)] Setting $p_1^{\hat Y} \coloneq p_\infty$, the curve $(p^{\hat Y}_\tau)_{\tau \in [0, 1]}$ is Lipschitz continuous in Wasserstein-2, i.e.\ it holds
       \begin{equation}
           W_2\left(p_{\tau_1}^{\hat Y}, p_{\tau_2}^{\hat Y}\right) \leq \left(\sqrt{\mathbb{E}[|Y_0|^2]} + \frac{c}{\theta}\right) |\tau_1 - \tau_2|  \quad \text{for all }\, \tau_1, \tau_2 \in [0, 1].
       \end{equation}
       If $\tau_1 = 1$ or $\tau_2 = 1$, we get the better constant $\mathbb{E}[|Y_0 - Y_\infty|^2]^\frac{1}{2}$ as in \eqref{Wasserstein-Bound}.
   \end{enumerate} 
\end{proposition}

\begin{remark}
    Notice the subtle difference between the above result and Proposition \ref{prop:lipschitz}: while the inhomogeneous Kac process $(Z_u)_{u \in [0,\infty)}$ is Lipschitz continuous in Wasserstein-2 \emph{and pathwise}, its compressed version $(\hat Z_\tau)_{\tau \in [0,1)}$ oscillates between $1$ and $-1$ \emph{infinitely often} as $\tau \to 1$ by the vanishing switching rate \eqref{specific-switching-OU-K}. Hence, although the compressed Kac process admits a Wasserstein-2 extension, there is no pathwise continuous extension to $\tau = 1$.
\end{remark}

Finally, we establish a boundedness result for the conditional velocity field of the rescheduled OU-K process $\tilde Y_{t} = Y_{g(t)} = \hat Y_{\tau(t)}$. It follows directly from the chain rule and Corollary \ref{cor:conditional-bound}.
\begin{corollary}\label{cor:rescheduled-velocity-bounded}
    Let $(\tilde Y_t)_{t \in [0, 1)}$ be the rescheduled OU-K process \eqref{1d-backbone} with rate $\theta > 0$ and speed $c > 0$, and differentiable time schedule $g:[0,1)\to [0,\infty)$ with $g(0)=0$ and $\lim_{t \uparrow 1} g(t) = \infty$. 
    Set $\tau \coloneqq \tau(t) \coloneqq 1 - e^{-\theta g(t)}, \,  t \in [0,1)$. 
    Then, its conditional velocity field \eqref{cond-velo-schedule} is given by
    \begin{equation}\label{rescheduled-velocity}
        v_t^{\tilde Y}(y \mid y_0) = \tau'(t) \left(- y_0 + \frac{c}{\theta}v_{\hat Z}(\tau(t), x) \right),
    \end{equation}
    where $x \coloneqq \frac{\theta}{c}(y - (1 - \tau(t))y_0)$, and $v_{\hat Z}$ is the velocity field \eqref{compressed-velocity-legendre}. Further, the conditional velocity admits the bound
    \begin{equation}\label{rescheduled-conditional-velo-bound} 
        |v^{\tilde Y}_t (y \mid y_0)|  \leq \tau'(t) \left(|y_0| + \frac{c}{\theta}\right) \quad \text{for all } t \in [0, 1).
    \end{equation}
\end{corollary}
Remarkably, the bound \eqref{rescheduled-conditional-velo-bound} becomes \emph{uniform} (in time and space), if we assume the data samples to be bounded, i.e.\ $|y_0| \le R$, and additionally assume 
\begin{equation}\label{weak-assumption}
    |\tau'(t)| = \theta |g'(t)| e^{-\theta g(t)} \le C < \infty \quad \text{for all } t \in [0,1).
\end{equation}
Note that the condition \eqref{weak-assumption} is \emph{weaker} than assuming that $|g'|$ is bounded. In particular, in our following generative OU-K model, we will use a time schedule \eqref{cosine-schedule} such that $|g'|$ is unbounded, but \eqref{weak-assumption} is still valid.
Hence, our OU-K model admits \emph{bounded} velocity fields, which render the usual truncations of diffusion models unnecessary.

\section{OU--K process in generative modeling}\label{sec:generative-modeling}

In the following, we describe our experimental setup and numerical results, demonstrating the performance of our OU-K generative model.

\subsection{OU--Kac generative model and FM/OU baseline}

\paragraph{The OU--K model.}
As motivated in Remark \ref{remark:time-schedule}, the backbone of our OU-K model is given by the rescheduled OU-K process in the form given in Proposition \ref{prop:transformation} by
\begin{equation}\label{backbone}
    \mathbf{\tilde Y}_t =  \mathbf{Y}_{g(t)} = e^{- \theta g(t)} \mathbf{X}_{0} + \frac{e^{-\theta g(t)}}{\theta} \, c \,  \mathbf{Z}_{e^{\theta g(t)} - 1}, \quad t \in [0,1),
\end{equation}
where $\mathbf{X}_0 \in \R^d$ is the random variable representing the target data distribution, $g: (0,1) \to \R$ is a suitable time-scheduling function, and $\mathbf{Z}_u = (Z_u^1, \dots, Z_u^d)$ consists of
$d$ i.i.d.\ unit-speed, inhomogeneous Kac processes $(Z_u^i)_u$ with switching rate given in \eqref{specific-switching-OU-K}. We choose $\mathbf{X}_0$ and $\mathbf{Z}_u$ to be independent.
Based on our experiments, we have found the commonly used \emph{cosine schedule}, see \cite{nichol2021improved}, to be well suited for our OU-K model. With a small fixed offset $s \approx 0.008$, it is defined by
\begin{equation}\label{cosine-schedule}
    g: [0,1) \to \R, \quad 
    g(t) = - 2 \ln \Big( \frac{\cos(\frac{\pi}{2} \frac{t+s}{1+s} )}{\cos(\frac{\pi}{2} \frac{s}{1+s} )} \Big).
\end{equation}
Here, we do not apply any time truncation to this schedule, see Appendix \ref{sec:implementation-details} for further details.  

\emph{Training:} 
In order to train a neural network to minimize the CFM loss \eqref{cfm_loss}, we sample $t \sim {\rm Unif}([0,1])$, $x_0 \sim \mathbf{X}_0$, and $z \sim \mathbf{Z}_{e^{\theta g(t)}-1}$, which gives us $y \sim {\rm Law} (\mathbf{Y}_{g(t)} \mid \mathbf{Y}_0 = x_0)$ via \eqref{backbone}.
For sampling $z \sim \mathbf{Z}_{e^{\theta g(t)}-1}$, we apply the explicitly given law \eqref{legendre-formula-total} derived in Theorem \ref{thm:main-polynomial-law}. The conditional velocity field $v^\mathbf{\tilde Y}_t(y \mid x_0)$ in the CFM loss \eqref{cfm_loss} is then computed by combining \eqref{Y-cond-velo}, \eqref{cond-velo-multi-d} and \eqref{cond-velo-schedule}.

\emph{Inference:}
Once the neural network $v_t^\theta$ is trained, we obtain new samples from the target distribution $p_0$ by integrating $-v_{1-t}^\theta$ in \eqref{eq:flow_ode}, and applying $p_0 = \varphi(1,\cdot)_\sharp p_\infty$, where $p_\infty$ is the invariant distribution from Theorem \ref{thm:invariantkacou}.

\paragraph{The standard flow matching (FM) baseline.}
As a baseline for our experiments, we use the well-known linear interpolation path given by
\begin{equation}
    \mathbf{X}_t^{\rm FM} \coloneqq (1-t) \mathbf{X}_0 + t \, \mathbf{\Xi}, \quad t \in [0,1],
\end{equation}
with a centered Gaussian $\mathbf{\Xi} \sim \mathcal{N}(0, I_d)$ independent of $\mathbf{X}_0$, with latent variance $\sigma_{\rm FM}^2 = 1$.
The model is also known as \emph{independent coupling flow matching}, for details, we refer to \cite{lipman2023flow, WS2025}.

\paragraph{The OU (diffusion) baseline.}
We also compare our OU-K model against an FM model, where the backbone is given by the rescheduled Ornstein-Uhlenbeck process
\begin{equation}
       \mathbf{X}_{g(t)}^{\rm OU} \overset{d}{=} 
       e^{-\theta g(t)} \mathbf{X}_0 + \sqrt{\frac{1 - e^{-2 \theta g(t)}}{2\theta}} \, \sigma \, \mathbf{\Xi},
\end{equation}
with a standard Gaussian $\mathbf{\Xi} \sim \mathcal{N}(0, I_d)$. We choose (a truncated version of) the cosine schedule $g(t)$ in \eqref{cosine-schedule}, the rate $\theta = 1$, and a volatility $\sigma$ such that for the terminal OU variance it holds
$\sigma_{{\rm OU}}^2 = \frac{\sigma^2}{2 \theta} \in \{\frac{1}{100}, 1\}$, see Appendix \ref{sec:implementation-details}.
The conditional velocity field is here given by \eqref{OU-cond-velocity}.

\paragraph{The Kac baseline.}
For completeness, we also include an FM model employing the mean-reverting Kac process
\begin{equation}
       \mathbf{M}_{t}^{\rm Kac} =
       (1-t) \mathbf{X}_0 + \mathbf{K}_{t^2},
\end{equation}
where the components of $\mathbf{K}_t$ are 1D homogeneous Kac processes with switching rate $a$ and speed $c$,
see \cite{duong2026telegraphersgenerativemodelkac}. Here, we choose the damping $a \in \{50, 100\}$ and adjust the Kac speed $c$ such that the latent variance $\frac{c^2}{a}
(1-\frac{1-e^{-2a}}{2a})$ matches the terminal FM variance $\sigma_{\rm FM}^2$.

\paragraph{Choosing the OU-K parameters.}
Our OU-K backbone $\mathbf{Y}_{g(t)}$ given in \eqref{backbone} is determined by the mean-reversion rate $\theta > 0$, speed $c>0$ and switching rate $a > 0$. In our experiments, we choose $\theta = 1$, and a suitable integer value $a \in \N$, such that the assumption 
$k = \frac{a}{\theta} = a \in \N$ of Theorem \ref{cor:velo-formula} is fulfilled. The speed $c > 0$ is then adjusted such that the latent variance $\sigma_{{\rm OU-K}}^2 = \frac{c^2}{\theta(\theta + 2a)}$, see Theorem \ref{thm:invariantkacou}, matches the OU baseline latent variance $\sigma_{\rm OU}^2$.

\paragraph{Choosing the terminal time $T$.}
Notice that with the choice $\theta = 1$, we have
\begin{equation}\label{transformed-time}
  \tau(t) \coloneq 1- e^{-g(t)} = 1- \frac{\cos^{2}(\frac{\pi}{2} \frac{t+s}{1+s})}{\cos^{2}(\frac{\pi}{2} \frac{s}{1+s})} \quad ~ \text{and} \quad ~ 
  \tau'(t) = \frac{\pi}{2(1+s)\cos^{2}(\frac{\pi}{2} \frac{s}{1+s})} \sin \left(\pi \frac{t+s}{1+s} \right).
\end{equation}
Since $\sup_t |\tau'(t)| < \infty$, the conditional velocity field $v^\mathbf{\tilde Y}_t(y \mid y_0)$ stays uniformly bounded for all $t \in [0,1)$ by \eqref{rescheduled-conditional-velo-bound}, given that the data $|y_0| \le R$ is bounded. 
Therefore, we have a theoretical boundedness guarantee for \emph{all} time points, and the time truncations near data and latent usually performed in diffusion models, see e.g.\ \cite{soft2021, nichol2021improved}, \emph{are not necessary} in our OU-K model. Hence, we can choose the initial time $T_0 = 0$ and terminal time $T_1=1$ without restrictions.
\\[1ex]
Further implementation details are given in Appendix \ref{sec:implementation-details}.

\subsection{The sharpness experiment on singular data sets}\label{sec:sharpness}

Diffusion models are known to blur out data sets with sharp, singular modes, see \cite{duong2026telegraphersgenerativemodelkac}. With this experiment, we aim to demonstrate that our OU-Kac model is able to generate such data sets with more precision and fewer training steps.

We consider the $d$-dimensional toy data sets
\begin{itemize}
    \item Cross-polytope: $D_1 = \{\pm e_1, \dots, \pm e_d\}$ with the unit vectors $e_1, \dots, e_d \in \R^d$.
    \item Rotated cross-polytope: $D_2 = U D_1$, i.e.\ $D_1$ rotated by some orthogonal matrix $U$.
    \item Simplex: $D_3 = \{x_0, \ldots, x_d\}$, a regular simplex centered at the origin.
    \item Unit sphere: $D_4 = \{x \in \R^d: \|x\|_2 = 1\}$.
\end{itemize}
In fact, we carry out the experiments on rescaled versions of $D_i$, see Appendix \ref{sec:implementation-details} for details.
For each data set, we consider a target distribution uniformly distributed on the points in $D_i$, $i = 1,2,3,4$.
We fix the dimension at $d = 128$ and study the generation performance of our OU-K model, the standard FM, OU and Kac baselines on these toy sets. For our OU-K model and the Kac baseline, we choose the damping parameter $a=50$, and adjust the other parameters as described in Appendix \ref{sec:implementation-details}.
We implemented our own metric for discrete data, the \emph{balanced AUC sharpness profile} \eqref{balanced_sharpness_AUC}, measuring the sharpness of the generated points, as well as the ability to uniformly cover the data sets; for the details, see Appendix \ref{sec:metric-discrete}. 

\begin{figure}[H]
    \centering
    \includegraphics[width = 1.0\textwidth]{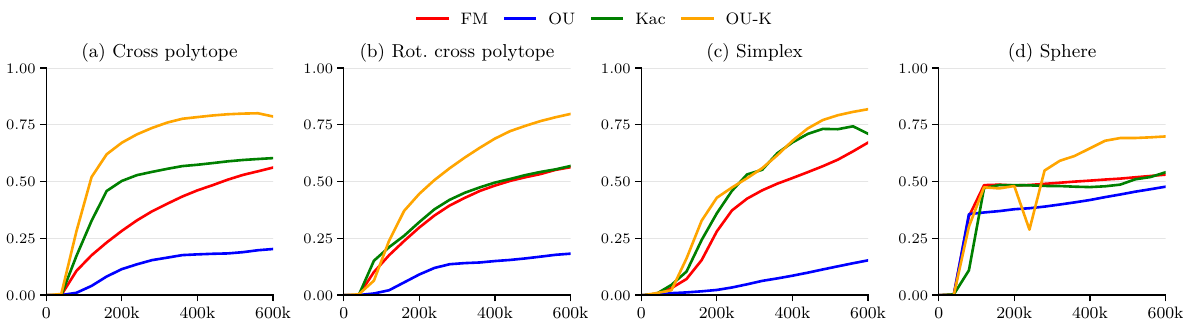}    
     \caption{The balanced AUC sharpness metric of the data sets $D_i$, $i=1,2,3,4$, in dimension $d=128$, as a function of the number of training steps, comparing the OU-K model with the FM, OU and Kac baselines.}
    \label{fig:four-sets-d128}
\end{figure}

\paragraph{Results.}
In Figure \ref{fig:four-sets-d128}, we report the evolution of our sharpness metric \eqref{balanced_sharpness_AUC} over the course of training. Overall, the OU-K model achieves a \emph{substantial improvement in sample sharpness} over the flow matching and diffusion baselines, reaching sharpness scores of $> 75 \, \%$ for $D_1,D_2,D_3$, and $> 70 \, \%$ for $D_4$ at 600k training steps. Figure \ref{fig:OU-K_vs_OU_intro} gives a visualization of the generated point clouds.
\\[1ex]
\textbf{Comparing OU-K with FM:}
\begin{itemize}
    \item Cross-polytope $D_1$: OU-K substantially outperforms the FM baseline in terms of sharpness, with an improvement of $\approx 40 \, \%pt.$ at 200k steps, and $\approx 30 \, \%pt.$ after 520k training steps.
    \item Rotated polytope $D_2$: the OU-K model achieves a significant sharpness improvement of $\approx 24 \, \%pt.$ over the FM baseline after 600k training steps.
    \item Simplex $D_3$: we observe a smaller performance increase of $\approx 15-20 \, \%pt.$ after 440k steps over the FM baseline, partly attributable to the latter's stronger performance on this dataset.
    \item Sphere $D_4$: despite a rather irregular initial phase, OU-K eventually achieves a sharpness improvement of $\approx 17 \, \%pt.$ over the FM baseline after 440k steps.
\end{itemize}
Moreover, the \textbf{OU-K model severely and consistently outperforms the OU baseline}, with an improvement of $> 50 \, \%pt.$ on $D_1,D_2,D_3$, and $\approx 25 \, \%pt.$ on $D_4$.
This observation is consistent with the uniform boundedness of the (conditional) OU-K velocity field established in Corollary \ref{cor:rescheduled-velocity-bounded}, in contrast to the blow-up of the OU velocity field, cf.\ \eqref{OU-cond-velocity} and Remark \ref{rem:difference-OU-K-OU}.

\begin{remark}
    We note that the substantial improvement in sample sharpness achieved by the OU-K model on the cross-polytope $D_1$ cannot be attributed solely to its componentwise construction. First, the results on the rotated polytope $D_2$ show that the OU-K model's advantage persists under a rotation of the data, providing a useful sanity check. Second, the flow matching baseline \emph{implicitly} incorporates a componentwise construction as well, see \cite[Appendix B.1]{CKDS2026}. 
    A possible structural distinction is rather that the OU-K model is \emph{non-radial}, in contrast to the FM baseline. However, the OU-K model also exhibits improved sample sharpness on the sphere $D_4$, suggesting that its gains may instead stem from the greater regularity of the OU-K model itself, see Corollary \ref{cor:rescheduled-velocity-bounded}.
\end{remark}

\subsection{CIFAR-10 image generation} \label{sec:straight-trajectories}

As a proof of concept, we compare our OU-K model against the flow matching, OU and Kac 
baselines on CIFAR-10 image generation. For the flow matching and Kac baselines, we fix the latent variance $\sigma^2_{\rm FM} = 1$, while for the OU and OU-K model, we have empirically found a small latent variance $\sigma^2_{\rm OU} = \frac{1}{100}$ to be preferable. 
For the Kac and OU-K model, we 
choose the damping parameter $a = 100$, and adjust the Kac speed $c$ accordingly to obtain the desired latent variance. 
We parameterize the velocity field in the CFM loss \eqref{cfm_loss} using a U-Net architecture and simulate the resulting flow ODE \eqref{eq:flow_ode} with 
standard 100- and 200-step Euler integration. 
The details are provided in Appendix \ref{sec:implementation-details}.
In Table \ref{tab:cifar-euler-variance0p01}, we report the Fréchet inception distance (FID) scores, and the corresponding generated images are shown in Figure \ref{fig:cifar-1}. 
We observe that our OU-K model achieves generation performance \emph{competitive} with standard generative flow models such as the FM and OU baselines.


\begin{table}[H]
\centering
\setlength{\tabcolsep}{8pt}
\begin{tabular}{c|cccc}
\toprule
FID $\downarrow$ & FM & OU & Kac & OU-K \\
\midrule
Euler 100 & 4.01 & 4.15 & 5.05 & 4.12 \\ 
Euler 200 & 3.59 & 3.50 & 4.10 & 3.48 \\
\bottomrule
\end{tabular}
\caption{CIFAR-10 FID scores evaluated at 50k samples. 
Lower score is better.}
\label{tab:cifar-euler-variance0p01}
\end{table}

\paragraph*{Acknowledgments.}
The work of R. Duong was funded by the Deutsche Forschungsgemeinschaft (DFG) under Germany’s Excellence Strategy – The Berlin Mathematics Research Center MATH+ (EXC-2046/1, EXC-2046/2, project ID: 390685689).
GS acknowledges funding by DFG, Project STE 571/23-1.

\bibliographystyle{abbrv}
\bibliography{Bibliography}

\appendix
\section{Additional Proofs}\label{sec:appendix}

\subsection{Proofs of Section \ref{sec:finite-OU-K}}\label{appendix-proofs-A1}

\begin{proof}[\textbf{Proof of Proposition~\ref{prop:kac-scaling}}]
Let $T > 0$ and $(K_t^{(n)})_{t \in [0, T]}$ be the symmetric and time-homogeneous Kac process \eqref{hom-Kac-process} driving the OU-K process $(Y_t^{(n)})_{t \in [0, T]}$ with speed $c_n$ and switching intensity $a_n$. Using the integration-by-parts formula for Riemann-Stieltjes integrals, it holds
\begin{equation}\label{ibp-RS}
    \int_0^t e^{- \theta(t - s)} \d K_s^{(n)} = K_t^{(n)} - \theta \int_0^t e^{- \theta (t - s)} K_s^{(n)} \d s,
\end{equation}
while the integration-by-parts formula for It{\^o} processes yields
\begin{equation}\label{ibp-ITO}
    \int_0^t e^{- \theta(t - s)} \d W_s = W_t - \theta \int_0^t e^{- \theta(t - s)} W_s \d s.
\end{equation}
We define the map
\begin{equation}
    \Phi: \R \times C([0, T]) \to C([0, T]),\ (x, f) \mapsto \Phi(x, f)
\end{equation}
with
\begin{equation}
    \Phi(x, f): [0, T] \to \R,\ t \mapsto e^{- \theta t}x + \mu(1 - e^{-\theta t})+ f(t) - \theta \int_0^t e^{- \theta(t - s)} f(s) \d s,
\end{equation}
such that with \eqref{ibp-RS} it holds
\begin{equation}\label{OU-K-as-cont-map-of-K}
    (Y_t^{(n)})_{t \in [0, T]} = \Phi(Y_0^{(n)}, (K_t^{(n)})_{t \in [0, T]}).
\end{equation}
Further, the map $\Phi$ is continuous since for any $(y_n, g_n)_{n \in \N} \in \R \times C([0, T])$ with $(y_n, g_n) \to (x, f)  \in \R \times C([0, T])$, it holds
\begin{align}
    \|\Phi(x, f) - \Phi(y_n, g_n)\|_{\infty} &= \sup_{t \in [0, T]} |\Phi(x, f)(t) - \Phi(y_n, g_n)(t)| \\ 
    &\leq \sup_{t \in [0, T]} e^{-\theta t} |x - y_n| + \|f - g_n\|_{\infty} + \theta \sup_{t \in [0, T]} \int_0^t e^{- \theta(t - s)} |f(s) - g_n(s)| \d s \\ 
    &\leq |x - y_n| + 2 \|f - g_n\|_{\infty} \to 0
\end{align}
as $n \to \infty$. With \eqref{OU-K-as-cont-map-of-K}, the Kac scaling $ \lim_{n \to \infty} (K_t^{(n)})_{t \in [0, T]} \overset{d}{=} (\sigma W_t)_{t \in [0, T]}$ \cite{KAC1974} and the continuous mapping theorem, we obtain
\begin{equation}
    \lim_{n \to \infty} (Y_t^{(n)})_{t \in [0, T]} = \lim_{n \to \infty} \Phi(Y_0, (K_t^{(n)})_{t \in [0, T]}) \overset{d}{=} \Phi(Y_0, (\sigma W_t)_{t \in [0, T]}) \eqcolon (X_t)_{t \in [0, T]}.
\end{equation}
By \eqref{ibp-ITO}, the process $(X_t)_{t \in [0, T]}$ satisfies
\begin{align}
    X_t &= e^{- \theta t}X_0 + \mu(1 - e^{-\theta t}) + \sigma W_t - \theta \sigma \int_0^t e^{- \theta (t - s)} W_s \d s \\ 
    &= e^{-\theta t} X_0 + \mu (1 - e^{-\theta t}) + \sigma \int_0^t e^{- \theta(t - s)} \d W_s,
\end{align}
which is exactly the solution of the SDE \eqref{Wiener-SDE}. 
\end{proof}

\begin{proof}[\textbf{Proof of Proposition~\ref{Z-Kac-is-MSP}}]
By definition \eqref{inhom-Kac-process-Zt}
we have that $(Z_t)_t$ is pathwise continuous.
We abbreviate $I_t \coloneqq I_t^Z$.
Next, denote the times of the Poisson events of $(N_t)_t$ by $T_k$, $k=1,2,\dots$, and without loss of generality, let $N_{T_k}$ be even, i.e.
\begin{align}
    I_{T_k} = \frac{1}{2} B_{\frac{1}{2}} (-1)^{N_{T_k}} + \frac{3}{2} = 
    \begin{cases}
        \, 1,  \quad \text{if } B_{\frac{1}{2}} = -1,\\ 
        \, 2,  \quad \text{if } B_{\frac{1}{2}} = 1.
    \end{cases}
\end{align}
Note that here, $F^{I_{T_k}}_Z(x) = B_{\frac{1}{2}} c_Z$.
The process $(I_t)_t$ is constant for $t \in [T_k, T_{k + 1})$, and by definition \eqref{inhom-Kac-process-Zt} the continuous evolution of $Z_t$ satisfies
\begin{equation}
    Z_t = Z_{T_k} + B_{\frac{1}{2}} c_Z (t - T_k), \quad t \in [T_k, T_{k + 1}),
\end{equation}
and therefore,
\begin{equation}
    \d Z_t = B_{\frac{1}{2}} c_Z \d t = F^{I_{T_k}}_Z(Z_t) \d t.
\end{equation}
The distribution of the inhomogeneous Poisson interarrival time $T_{k + 1} - T_k$ is
\begin{equation}
    \mathbb{P}(T_{k + 1} - T_k > t | Z_{T_k} = z, I_{T_k} = i) = \exp \left( -\int_{T_k}^{T_k + t} a_Z(s) \d s \right) 
    = \exp\left(- \int_0^t \lambda_Z^{i}(T_k + s, \phi_s^i(x)) \d s \right),
\end{equation}
where $\phi_s^i(z) \coloneqq z+ (-1)^i c_Z s$, since by definition, $\lambda^{i}_Z(t, x) = a_Z(t)$ for $i=1,2$.
At the next Poisson event $T_{k + 1}$, the term $N_{T_{k + 1}}$ becomes odd and it follows
\begin{equation}
    I_{T_{k + 1}} = 
    \begin{cases}
        \, 2,  \quad \text{if } B_{\frac{1}{2}} = -1,\\ 
        \, 1,  \quad \text{if } B_{\frac{1}{2}} = 1.
    \end{cases}
\end{equation}
Hence a jump in $N_t$ corresponds to a switch in $I_t$, 
and $I_{T_{k + 1}} \ne I_{T_{k}}$ almost surely.

Lastly, it holds for $i, j = 1, 2$ that
\begin{align}
    \mathbb{P}(I_{T_{k + 1}} = j \mid Z_{T_{k}} = z, I_{T_k} = i) 
    \,=\, \mathds 1_{\{j \ne i\}}
    \,=\, \frac{\lambda_Z^{ij}(T_{k + 1}, z)}{\lambda_Z^i(T_{k + 1}, z)},
\end{align}
since $\lambda^{i}_Z(t, x) = a_Z(t)$ and 
$\lambda_Z^{ij}(t, x) = a_Z(t) \mathds 1_{\{j \ne i\}}$. 
\end{proof}

\begin{proof}[\textbf{Proof of Proposition~\ref{FPE-MSP}}]
Let $\psi \in C_c^\infty((0, \infty) \times \R^d)$ and $k \in \{1, \ldots, n\}$. Define 
\begin{equation}\label{test-function-construction}
    \varphi: (0, \infty) \times \R^d \times \{1, \ldots, n\} \to \R,\ (t, x, i) \mapsto \mathds{1}_{\{i = k\}}\psi(t, x).
\end{equation}
Then $(\varphi, \mathcal{A} \varphi)$ with
\begin{equation}\label{extended-generator}
    \mathcal{A} \varphi(t, x, i) = \partial_t \varphi(t, x, i) + \langle F^i(x), \nabla_x \varphi(t, x, i)\rangle + \sum_{j = 1}^n \lambda^{ij}_t(x)(\varphi(t, x, j) - \varphi(t, x, i))
\end{equation}
is in the extended generator of $(t, X_t, I_t)_t$, see \cite[Theorem 5.5]{davis1984}, i.e. the process $(M_t^\varphi)_t$ with
\begin{equation}\label{extended-generator-martingale}
    M_t^\varphi \coloneq \varphi(t, X_t, I_t) - \varphi(0, X_0, I_0) - \int_0^t \mathcal{A} \varphi(s, X_s, I_s) \d s
\end{equation}
is a martingale with respect to the filtration $(\sigma(X_s, I_s \mid 0 \leq s \leq t))_t$. Since $\varphi$ is compactly supported in $(0, \infty) \times \R^d$, the first summand in \eqref{extended-generator-martingale} vanishes for sufficiently large $t \geq 0$, while the second one vanishes by compact support of $\varphi$. Further, as $\varphi$ is only supported on the mode $k$, \eqref{extended-generator} simplifies to
\begin{equation}\label{extended-generator-specific}
    \mathcal{A} \varphi(t, x, i) = \begin{cases}
        \partial_t \psi(t, x) + \langle F^k(x), \nabla_x \psi(t, x)\rangle - \sum_{j = 1}^n \lambda^{kj}_t(x)\psi(t, x) & \text{for } i = k,\\ 
        \lambda^{ik}_t(x) \psi(t, x) & \text{else}.
    \end{cases}
\end{equation}
Taking the expectation and $t \to \infty$ in \eqref{extended-generator-martingale} and using the representation \eqref{extended-generator-specific} yields 
\begin{align}
    0 = \int_0^\infty \mathbb{E}[\mathcal{A} \varphi(s, X_s, I_s)] \d s &= \int_0^\infty \sum_{i = 1}^n \int_{\R^d} \mathcal{A} \varphi(t, x, i) \d \rho_{i, t}(x) \\ 
    &= \int_0^\infty \int_{\R^d} \partial_t \psi(t, x) + \langle F^k(x), \nabla_x \psi(t, x) \rangle \d \rho_{k, t}(x)\d t \\ 
    &\phantom{=} + \sum_{i \neq k} \int_0^\infty \left(\int_{\R^d} \lambda^{ik}_t\psi(t, x) \d \rho_{i, t}(x) - \int_{\R^d} \lambda^{ki}_t \psi(t, x) \d \rho_{k, t}(x)\right) \d t.
\end{align}
This last term is the distributional formulation of the PDE
\begin{equation}
    \partial_t \rho_{k, t} = - \operatorname{div}(F^k \rho_{k, t}) + \sum_{i \neq k} (\lambda^{ik}_t \rho_{i, t} - \lambda^{ki}_t \rho_{k, t}),
\end{equation}
which finishes the proof. 
\end{proof}

\noindent
\textbf{Addition to the proof of Theorem \ref{prop:PDE-1} ( Distributional solutions).}
Since the solutions of the damped wave equation \eqref{damped-wave} admit singular parts, we need to understand them in a distributional sense. 

Given initial distributions $f,g \in \mathcal{D}'(\R)$, we say that a weakly continuous curve $p:[0,\infty) \to \mathcal{P}(\R), \, t \mapsto p_t = p(t,\cdot)$, solves 
\begin{align}\label{damped-wave-distributional-1}
    \begin{cases}
    \partial_{tt}p + 2 a(t) \partial_{t}p \,=\,  c^2 \partial_{xx} p,\\ 
    p(0, \cdot) = f,\\ 
\partial_{t}p(0, \cdot) = g,\end{cases}
\end{align}
in the sense of distributions iff it holds
\begin{align}\label{damped-wave-distributional-2}
    \int_0^\infty \int_{\R} (\partial_{tt}\varphi(t,x) - 2 a(t) \partial_t \varphi(t,x))  \d p_t(x) \d t
    = c^2 \int_0^\infty \int_{\R} \partial_{xx}\varphi(t,x)  \d p_t(x) \d t,
\end{align}
for all smooth test functions $\varphi \in C_c^\infty((0,\infty) \times \R)$, and the following initial conditions are satisfied: 
\begin{itemize}
    \item $p(0,\cdot) = f$ (as distributions).
    \item For any family of smoothing kernels $j_t^\varepsilon \in C_c^\infty( (t-\varepsilon, t+\varepsilon))$ such that $\lim_{\varepsilon \to 0} \int_0^\infty j_t^\varepsilon(s) h(s) \d s = h(t)$ for any continuous $h:(0,\infty) \to \R$, 
    it holds
    \begin{equation}
        \lim_{t \to 0} \lim_{\varepsilon \to 0} - \int_0^\infty \partial_s j_t^\varepsilon(s) \int_{\R} \psi(x) \d p_s(x) \d s
        = g(\psi)
    \end{equation}
    for all $\psi \in C_c^\infty(\R)$.
\end{itemize}

\paragraph{Completing the proof.}
The equation \eqref{damped-wave-distributional-2} follows directly by the formal computations made in the main-part proof of Theorem \ref{prop:PDE-1}, applied to the test function $\varphi$.

We only need to check the initial condition $\partial_{t}p(0, \cdot) = 0$ in the distributional sense: For this, fixing $\psi \in C_c^\infty(\R)$, we see that $j_t^\varepsilon(s) \psi(x) \in C_c^\infty((0,\infty) \times \R)$ for small enough $\varepsilon$. Therefore, by \eqref{KFE-mode} it holds
\begin{align}
   & \int_0^\infty \int_{\R} \partial_s j_t^\varepsilon(s) \psi(x) \d p_s(x) \d s\\
  =~ &\int_0^\infty \int_{\R} \partial_s j_t^\varepsilon(s) \psi(x) \d \rho_{1,s}(x) \d s
    + \int_0^\infty \int_{\R} \partial_s j_t^\varepsilon(s) \psi(x) \d \rho_{2,s}(x) \d s\\
  =~ & c\int_0^\infty \int_{\R} j_t^\varepsilon(s) \partial_x\psi(x) \d \rho_{1,s}(x) \d s
    - c \int_0^\infty \int_{\R} j_t^\varepsilon(s) \partial_x \psi(x) \d \rho_{2,s}(x) \d s\\  
  \longrightarrow ~ & c\int_{\R} \partial_x\psi(x) \d \rho_{1,t}(x)
    - c\int_{\R} \partial_x\psi(x) \d \rho_{2,t}(x)
\end{align}
for $\varepsilon \to 0$. Finally, the weak continuity of $t \mapsto \rho_{1,t}$, $t \mapsto \rho_{2,t}$, together with $\rho_{1,0} = \rho_{2,0} = \frac{1}{2} \delta_0$,  yields
\begin{align}
     \lim_{t \to 0} \lim_{\varepsilon \to 0} \int_0^\infty \int_{\R} \partial_s j_t^\varepsilon(s) \psi(x) \d p_s(x) \d s
     = \frac{c}{2} \partial_x\psi(0) 
    - \frac{c}{2} \partial_x\psi(0) = 0.
\end{align}

The weak continuity of $t \mapsto \rho_{i,t}$ can be seen as follows. Let $f \in C_b(\R)$ and fix a sequence $t_k \to t$ in $[0,\infty)$. Notice that $(Z_t)$ is almost surely Lipschitz continuous, and that for $\mathbb{P}$-a.e.\ $\omega$ it holds $I_{t_k}(\omega) = I_{t}(\omega)$ for $k$ large enough (depending on $\omega$)\footnote{$(I_t)_t$ is almost surely continuous, but has no continuous modification.}. 
Therefore, $f(Z_{t_k})\mathds{1}_{\{I_{t_k} = i \}} \to f(Z_{t})\mathds{1}_{\{I_{t} = i \}}$ $\mathbb{P}$-almost surely, and  
\begin{equation}
    \int_{\R} f(x) \d \rho_{i,t_k}(x) 
    = \int_\Omega f\big(Z_{t_k}(\omega)\big) \mathds{1}_{\{I_{t_k}(\omega) = i \}} \d \mathbb{P}(\omega)
    \to \int_\Omega f\big(Z_{t}(\omega)\big) \mathds{1}_{\{I_{t}(\omega) = i \}} \d \mathbb{P}(\omega)
    = \int_{\R} f(x) \d \rho_{i,t}(x),
\end{equation}
as $k \to \infty$, using Lebesgue's dominated convergence theorem.
\qed
\\[2ex]
\textbf{Addition to the proof of Theorem \ref{thm:main-polynomial-law} (Riemann-Volterra method).}
Let us consider the EPD equation \eqref{EPD} for general initial datums, and importantly, for \emph{positive} initial times $\tau_0 > 0$, that is, 
\begin{align}\label{EPD-Riemann-Volterra}
    \begin{cases}
        \partial_{\tau \tau} p + \frac{2 k}{\tau} \partial_{\tau}p ~=~  \partial_{xx} p, \quad ~ \tau > \tau_0 > 0,\\ 
        p(\tau_0, x) = \phi(x),\\ 
        \partial_{\tau}p(\tau_0, x) = \psi(x).\end{cases}
    \end{align}
We now derive the solution formula \eqref{EPD-solution} via the Riemann method. For further details, we refer to \cite[Chapter 5]{Copson_1975}, and also to the overview paper \cite{Zeitsch2018}.
\\[1ex]    
\emph{Step 1: Transformation of \eqref{EPD-Riemann-Volterra} to Canonical Form.}
Set $q(\tau,x) \coloneqq \tau^{k} p(\tau, x)$. A simple calculation shows that $p$ solves \eqref{EPD-Riemann-Volterra} iff $q$ solves the canonical form 
\begin{equation}\label{EPD-canonical}
    \partial_{\tau \tau} q - \frac{k(k-1)}{\tau^2} q ~=~  \partial_{xx} q, \quad ~ \tau > \tau_0,
\end{equation}
with the initial conditions (note $\tau_0 > 0$)
\begin{equation}\label{initial-canonical}
    q(\tau_0, x) = \tau_0^k \phi(x), \quad
    \partial_{\tau} q(\tau_0, x) = k \tau_0^{k-1} \phi(x) + \tau_0^k \psi(x).
\end{equation}
\emph{Step 2: The Riemann function.}
Define the operator corresponding to the canonical form \eqref{EPD-canonical} as
\begin{equation}
    L(q) \coloneqq \partial_{\tau \tau} q - \frac{k(k-1)}{\tau^2} q - \partial_{xx} q.
\end{equation}
The operator $L$ is self-adjoint, i.e.\ $L = L^*$, see \cite[Chapter 5.1, page 78]{Copson_1975}. The Riemann function $R(\tau, x; \tau_0, \xi)$ is given by the characteristic boundary-value problem
\begin{align}\label{Riemann-characteristic-problem}
    &L^*_{(\tau_0, \xi)}(R) \,=\, \partial_{\tau_0 \tau_0} R - \frac{k(k-1)}{\tau_0^2} R - \partial_{\xi \xi} R ~=~ 0,\\
    &R(\tau, x; \tau_0, \xi) \,=\, 1 ~~ \text{along the characteristic lines}~ \xi- x = \pm (\tau_0 - \tau).
\end{align}
For the problem \eqref{Riemann-characteristic-problem}, it can be given in terms of Gauss's hypergeometric function ${_2 F_1}$ via 
\begin{equation}\label{Gauss-hypergeometric}
    R(\tau, x; \tau_0, \xi) = {_2 F_1} (k, 1-k; 1; \sigma),
\end{equation}
where $\sigma \coloneqq \frac{(\tau - \tau_0)^2 - (x - \xi)^2}{4 \tau \tau_0}$, see \cite{Zeitsch2018}. 
Since in our case, $k \in \N$, this can also be written via the Legendre polynomial
\begin{equation}
    R(\tau, x; \tau_0, \xi) = L_{k-1}(1-2\sigma) = L_{k-1}\Big( \frac{\tau^2 + \tau_0^2 - (x-\xi)^2 }{2 \tau \tau_0}\Big).
\end{equation}
Note that indeed, we have $R (\tau, x; \tau_0, x \pm (\tau-\tau_0)) = L_{k-1}(1) = 1$.
\\[1ex] 
\emph{Step 3: The solution formula.}
Let 
\begin{align}
    H & \coloneqq H(\tau, x; \tau_0, \xi) = R \, \partial_\tau q - q \, \partial_{\tau_0} R\\
    & \coloneqq R(\tau, x; \tau_0, \xi) \, \partial_\tau q(\tau_0, \xi) - q(\tau_0, \xi) \, \partial_{\tau_0} R(\tau, x; \tau_0, \xi).
\end{align}
By \cite[Chapter 5.3, page 81]{Copson_1975}, it holds the solution formula
\begin{align}
    q(\tau, x) 
    &= \frac{1}{2} \Big( q (\tau_0, A) \, R (\tau, x; \tau_0, A) + q (\tau_0, B) \, R (\tau, x; \tau_0, B) \Big)  + \frac{1}{2} \int_{A}^{B}
    H   \d \xi\\
    &= \frac{1}{2} \Big( q (\tau_0, x - (\tau-\tau_0))  + q (\tau_0, x + (\tau-\tau_0)) \Big) \\
    &+ \frac{1}{2} \int_{x - (\tau-\tau_0)}^{x + (\tau-\tau_0)}
        R(\tau,x; \tau_0, \xi) \, \partial_\tau q(\tau_0, \xi) - q(\tau_0, \xi) \, \partial_{\tau_0} R(\tau,x; \tau_0, \xi) \d \xi,
\end{align}
where $A \coloneqq x - (\tau-\tau_0), \, B \coloneqq x + (\tau-\tau_0)$.
Dividing by $\tau^k$, inserting the initial conditions \eqref{initial-canonical}, and noticing that $p = \tau^{-k} q$,  yields the solution formula \eqref{EPD-solution}.
\qed

\begin{proof}[\textbf{Proof of Corollary~\ref{cor:law_OU-K}}]
    By Proposition \ref{prop:transformation}, we have 
    \begin{equation}
        Y_{t} = e^{- \theta t} X_{0} + \frac{e^{-\theta t}}{\theta} \, c \,  Z_{e^{\theta t} - 1}.
    \end{equation}
    We assume that the law of $X_0$ is absolutely continuous with continuous density $p_0$. Then, the law of $\tilde X_t \coloneqq e^{- \theta t} X_{0}$ reads as
    \begin{equation}
        p_{\tilde X}(t,x) = e^{\theta t} p_0(e^{\theta t} x).
    \end{equation}    
    By Theorem \ref{thm:main-polynomial-law}, the law of $Z_{e^{\theta t} - 1}$ is given by
    \begin{equation}
        p_Z(u, x) = \frac{1}{2(1+u)^k}(\delta(x-u) + \delta(x+u)) + 
        P_Z^{(k)}(u,x) \mathds{1}_{\{|x| < u\}},
    \end{equation}
    where $u \coloneqq e^{\theta t} - 1$, and $P_Z^{(k)}$ is given by the Legendre formula \eqref{legendre-formula-abs}. 
    Then, the law of $\tilde Z_t \coloneqq \frac{e^{-\theta t}}{\theta} \, c \,  Z_{e^{\theta t} - 1}$ is given by
    \begin{equation}
        p_{\tilde Z}(t,x) = \frac{\theta e^{\theta t}}{c} p_Z \left( u, \frac{\theta e^{\theta t}}{c} x \right).
    \end{equation}
    Finally, since $X_0$ and $Z_u$ are independent, the law of $Y_t$ is given by the convolution
    \begin{align}
        p_Y(t,x) &= \int_{-\infty}^\infty p_{\tilde X}(t,y) p_{\tilde Z}(t,x-y) \d y\\
        &= \frac{1}{2(1+u)^k} \Big\{ p_{\tilde X} \left(t,x - \frac{c(1-e^{-\theta t})}{\theta} \right) + p_{\tilde X} \left(t,x + \frac{c(1-e^{-\theta t})}{\theta} \right) \Big\} \\
        &+ \int_{x - \frac{c(1-e^{-\theta t})}{\theta}}^{x + \frac{c(1-e^{-\theta t})}{\theta}} p_{\tilde X}(t,y) \frac{\theta e^{\theta t}}{c} P_{Z}^{(k)}(u, \frac{\theta e^{\theta t}}{c}(x-y)) \d y\\
        &= \frac{e^{\theta t}}{2(1+u)^k} \Big\{ p_{0} \left(e^{\theta t} (x - \frac{c(1-e^{-\theta t})}{\theta} ) \right) + p_{0} \left(e^{\theta t} (x + \frac{c(1-e^{-\theta t})}{\theta}) \right) \Big\} \\
        &+ \frac{\theta e^{2 \theta t}}{c} \int_{x - \frac{c(1-e^{-\theta t})}{\theta}}^{x + \frac{c(1-e^{-\theta t})}{\theta}} p_{0}(e^{\theta t}y) \, P_{Z}^{(k)}(u, \frac{\theta e^{\theta t}}{c}(x-y)) \d y.
    \end{align}  
\end{proof}

\subsection{Proofs of Section \ref{sec:velocity-kac}}\label{appendix-proofs-A2}

\textbf{Addition to the proof of Proposition \ref{MSP-velo}.}
Concerning the absolute continuity of $(p_t)_t$ in $\mathcal{P}_2(\R^d)$, let us first show that $(p_t)_t$
is narrowly continuous, and it holds $p_t \in \mathcal{P}_2(\R^d)$ for all $t \geq 0$.

For the narrow continuity, notice that for all $f \in C_b(\R^d)$ and $(t_k)_{k \in \N}$ with $\lim_{k \to \infty} t_k = t \geq 0$, it holds
\begin{equation}
    \lim_{k \to \infty} \int_{\R^d} f(x) \d p_{t_k}(x) = \lim_{k \to \infty} \mathbb{E}[f(X_{t_k})] = \mathbb{E}[f(X_t)] = \int_{\R^d} f(x) \d p_t(x)
\end{equation}
by dominated convergence, and since the paths $t \mapsto X_t$ are continuous by construction.

Next, define $M \coloneq \max_{i \in \{1, \ldots, n\}} \|F^i(0)\|$ and $L \coloneq \max_{i \in \{1, \ldots, n\}} L^i$, where $L^i$ is the Lipschitz constant in \eqref{ass:Lipschitz}. Since
\begin{equation}\label{MSP-representation}
    X_t = X_0 + \int_0^t F^{I_s}(X_s)\d s,
\end{equation}
Bochner's inequality yields
\begin{equation}
    \|X_t\| \leq \|X_0\| + tM + L\int_0^t \|X_s\| \d s.
\end{equation}
With Grönwall's inequality, it follows
\begin{align}\label{expected-square-norm}
    \mathbb{E}[\|X_t\|^2] \leq 2e^{2Lt}(\mathbb{E}[\|X_0\|^2] + M^2t^2) < \infty,
\end{align}
and therefore $p_t \in \mathcal{P}_2(\R^d)$ for all $t \geq 0$. 

Lastly, we need to show $\|v_t\|_{L^2(\R^d, p_t)} \in L^2([0, T])$ for arbitrary $0 < T < \infty$. Let $t \geq 0$.  With conditional Jensen's inequality and the tower property of the conditional expectation it follows
\begin{align}\label{L2-norm-of-velo}
    \|v_t\|_{L^2(\R^d, p_t)}^2 &= \int_{\R^d} \|\mathbb{E}[F^{I_t}(X_t) \mid X_t = x]\|^2 \d p_t(x) \\ 
    &\leq \int_{\R^d} \mathbb{E}[\|F^{I_t}(X_t)\|^2 \mid X_t = x] \d p_t(x) \\ 
    &= \mathbb{E}[\|F^{I_t}(X_t)\|^2].
\end{align}
With $M$ and $L$ as above, Cauchy-Schwarz and Young's inequalities yield
\begin{align}\label{velo-bound}
    \mathbb{E}[\|F^{I_t}(X_t)\|^2] &= \mathbb{E}[\|F^{I_t}(X_t) - F^{I_t}(0)\|^2 + 2 \langle F^{I_t}(X_t) - F^{I_t}(0), F^{I_t}(0)\rangle + \|F^{I_t}(0)\|^2] \\ 
    &\leq L^2 \mathbb{E}[\|X_t\|^2] + 2 LM\mathbb{E}[\|X_t\|] + M^2 \\ 
    &\leq 2L^2 \mathbb{E}[\|X_t\|^2] + 2M^2.
\end{align}
Using \eqref{expected-square-norm} we arrive at
\begin{equation}\label{L2-norm-of-velo-bound}
\|v_t\|_{L^2(\R^d, p_t)}^2 \leq 4L^2e^{2Lt}(\mathbb{E}[\|X_0\|^2] + M^2t^2) + 2M^2,
\end{equation}
thus $\| v_t\|_{L^2(\R^d, p_t)} \in L^2([0, T])$. Together with the narrow continuity of $(p_t)_t$ and \eqref{eq:MSP-velo-CE}, it follows $(p_t)_t \in AC^2([0,T];\allowbreak \mathcal P_2(\R^d))$, see \cite[Theorem 1.1.2 \& 8.3.1]{BookAmGiSa05}. \qed

\begin{proof}[\textbf{Proof of Lemma~\ref{lem:legendre-quotient-estimate}}]
    1. By telescoping the recurrence formula \cite[Formula 28.23]{Spiegel2009}
    \begin{equation}
       L_{j+1}'(x) - L_{j-1}'(x) = (2j+1) L_j(x),
    \end{equation}
    and using $L_0' \equiv 0, \, L_1' \equiv 1 \equiv L_0$, we obtain for $k \ge 1$ the expansion formula\footnote{Technically, the expansion formula is still valid for $k=0$, yielding an empty sum in agreement with $L_0' \equiv 0$.}
    \begin{equation}
        L_k'(x) = \sum_{\substack{0\le j \le k-1\\j+k \text{ odd}}} (2j+1) L_j(x).
    \end{equation}
    Now, since $L_j(z) > 0$ for all $z \ge 1$ and $j\in \N_0$, it follows $L_k'(z) > 0$ for all $z \ge 1$, $k \ge 1$.
    \\[1ex]
    2. Hence, in order to prove \eqref{legendre-quotient}, we need to show that 
    \begin{equation}
       f(z) \coloneqq  L_{k}'(z) - L_{k - 1}'(z) > 0 \quad \text{for all } z \ge 1.
    \end{equation}
    With the recurrence relation \cite[Formula 28.22]{Spiegel2009}
    \begin{equation}
        L_{k-1}'(x) =  x L_{k}'(x) - k L_k(x),
    \end{equation}
    the claim is equivalent to
    \begin{equation}
        f(z) = L_k'(z) - z L_k'(z) + k L_k(z) > 0 \quad \text{for all } z \ge 1.
    \end{equation}
    Towards this end, first note that $f(1) = k L_k(1) = k > 0$. Next, the derivative of $f$ computes as
    \begin{equation}
        f'(z) = L_k''(z) - z L_k''(z) - L_k'(z) + k L_k'(z) 
        = (1-z) L_k''(z) + (k-1)L_k'(z).
    \end{equation}
    The Legendre differential equation \cite[Formula 28.1]{Spiegel2009}
    \begin{equation}
        (1-x^2) L_k''(x) = 2x L_k'(x) - k(k+1) L_k(x)
    \end{equation}
    yields via $1-x^2 = (1-x)(1+x)$ the identity
    \begin{equation}
        (1-x) L_k''(x) = \frac{2x L_k'(x) - k(k+1) L_k(x)}{1+x}.
    \end{equation}
    Plugging this into the derivative of $f$ results in
    \begin{equation}
        f'(z) = \frac{2z L_k'(z) - k(k+1) L_k(z)}{1+z} + (k-1)L_k'(z)
        = \frac{ [(1+z)(k-1) + 2z]  L_k'(z) - k(k+1) L_k(z)}{1+z}.
    \end{equation}
    Since $(1+z)(k-1) + 2z = z(k+1) + (k-1)$, this equals
    \begin{equation}
        f'(z) 
        = \frac{ z(k+1)L_k'(z) + (k-1)L_k'(z) - k(k+1) L_k(z)}{1+z}
        = \frac{ (k+1)(zL_k'(z) - k L_k(z)) + (k-1)L_k'(z)}{1+z}.
    \end{equation}
    Applying again the recurrence equation \cite[Formula 28.22]{Spiegel2009} yields
    \begin{equation}
        f'(z) = \frac{ (k+1) L_{k-1}'(z) + (k-1)L_k'(z)}{1+z}.
    \end{equation}
    Since by the first part, $L_k'(z), L_{k-1}'(z) \ge 0$ for all $z \ge 1$, the function $f$ is monotonically increasing on $[1, \infty)$ starting at $f(1) = k > 0$, and hence positive for all $z \ge 1$, which proves the claim \eqref{legendre-quotient}.
\end{proof}
Indeed, the above proof shows that the difference $f(z) = L_{k}'(z) - L_{k - 1}'(z)$, $z \in [1,\infty)$, is constantly $\equiv 1$, if $k=1$, and strictly increasing, if $k \ge 2$.
\\[2ex]
\textbf{Addition to Remark \ref{rem:difference-OU-K-OU} (Blow-up of OU velocity).}
With the schedule $g : [0,1) \to [0,\infty)$ and $\tau(t) \coloneqq 1 - e^{-\theta g(t)}$, the rescheduled version of the OU process \eqref{def:OU} (with mean $\mu=0$) reads
\begin{equation}
    \tilde X_t \coloneqq X_{g(t)} = \hat X_{\tau(t)},
\end{equation}
where the compressed process $(\hat X_\tau)_\tau$ is given in Remark \ref{rem:difference-OU-K-OU}.
Given that $g$ is differentiable, the conditional velocity field is then given by
\begin{equation}
    v_t^{\tilde X}(x\mid x_0) = \tau'(t) v_{\tau(t)}^{\hat X}(x\mid x_0)
    = \frac{\theta g'(t) e^{- \theta g(t)}}{1-e^{- 2 \theta g(t)}}
    \left(e^{- \theta g(t)} x - x_0\right).
\end{equation}
Let us study the term
\begin{equation}
    f(t) \coloneqq \frac{g'(t)}{1-e^{- 2 \theta g(t)}} = \frac{g'(t)}{g(t)} \cdot \frac{g(t)}{1-e^{- 2 \theta g(t)}}
\end{equation}
as $t \to 0$. For this, we assume that $g(0) = 0$, $g(t) > 0$ for all $t>0$, and $g$ is differentiable. Then, by L'Hospital's rule, 
\begin{equation}
    \lim_{t \to 0} \frac{g(t)}{1-e^{- 2 \theta g(t)}} =  \lim_{x \to 0} \frac{x}{1-e^{- 2 \theta x}} = \lim_{x \to 0} \frac{1}{2 \theta e^{- 2 \theta x}} = \frac{1}{2\theta}.
\end{equation}
Now, towards a contradiction, assume that for some $\varepsilon > 0$, it holds the bound
\begin{equation}
    \sup_{0 < t < \varepsilon} \left | \frac{g'(t)}{g(t)} \right | = C < \infty.
\end{equation}
Then, it follows
\begin{equation}
    g'(t) \le |g'(t)| \le C |g(t)| = C g(t)
\end{equation}
for all $0 < t < \varepsilon$. By Grönwall's inequality, we infer for all $0 < t < \varepsilon$ that
\begin{equation}
    g(t) \le g(0) e^{C t} = 0,
\end{equation}
contradicting the positivity assumption on $g$ for $t > 0$. This shows $\limsup_{t\to 0} |f(t)| = \infty$.

Hence, given that $x \ne x_0$, the rescheduled OU velocity field admits the blow-up
\begin{equation}
    \limsup_{t\to 0} |v_t^{\tilde X}(x\mid x_0)| = \limsup_{t\to 0} \left | \frac{\theta g'(t) e^{- \theta g(t)}}{1-e^{- 2 \theta g(t)}}
    \left(e^{- \theta g(t)} x - x_0\right) \right | 
    =\infty.
\end{equation}
In particular, this holds true for schedules of the form $g(t) = \int_0^t \beta(s) \d s$, where $\beta$ is continuous (or integrable) with $\beta(t) > 0$ for $t > 0$ and (possibly) $\beta(0) = 0$.

\begin{proof}[\textbf{Proof of Proposition~\ref{prop:compressed-W2-extension}}]
We use the shorthand $\left( a \right)_\# \mu \coloneq (T_a)_\# \mu$, where $T_a: \R \to \R,\ x \mapsto ax$, for $a \in \R$. 
\\
(i) We begin by showing that $(p_\tau^{\hat Z})_{\tau \in [0, 1)}$ is a 1-Lipschitz curve in $W_2$, from which the existence of a unique continuous extension to $\tau=1$ follows. From \eqref{compressed-Z-velo-bound} and \cite[Theorem 8.3.1]{AGS2008}, we may bound the metric derivative $|(p^{\hat Z})'|(\tau)$ of $(p_\tau^{\hat Z})_{\tau \in [0, 1)}$ by 
\begin{equation}
|(p^{\hat Z})'|(\tau) \leq \| v_{\hat Z}(\tau, \cdot)\|_{L^2(p_\tau^{\hat Z})} \leq 1 \qquad \text{for almost all }\tau \in [0, 1).
\end{equation}
Hence, using \cite[Theorem 1.1.2]{AGS2008}, we obtain the Lipschitz estimate
\begin{equation}\label{compressed-Z-Wasserstein}
    W_2(p_{\tau_1}^{\hat Z}, p_{\tau_2}^{\hat Z}) \leq \int_{\tau_1}^{\tau_2} |(p^{\hat Z})'|(r) \d r \leq \int_{\tau_1}^{\tau_2} \|v_{\hat Z}(r, \cdot)\|_{L^2(p_r^{\hat Z})} \d r \leq {\tau_2} - {\tau_1} 
\end{equation}
for all $0 \leq \tau_1 \leq \tau_2 < 1$.
By the completeness of $\mathcal{P}_2(\R)$, we can extend the Lipschitz curve $(p_\tau^{\hat Z})_{\tau \in [0, 1)}$ to the endpoint 
\begin{equation}
    p_1^{\hat Z} \coloneqq W_2 \text{ -} \lim_{\tau \to 1} p_{\tau}^{\hat Z}.
\end{equation}
To identify this endpoint, we rewrite \eqref{compressed-processes} to obtain
\begin{equation}
    \hat Z_\tau = \frac{\theta}{c}\left(\hat Y_\tau - (1 - \tau)Y_0\right).
\end{equation}
As $\hat Y_\tau \to Y_\infty$ in distribution and $(1 - \tau)Y_0 \to 0$ in probability, Slutsky's theorem yields $\hat Z_\tau \to \frac{\theta}{c} Y_\infty \sim \left( \frac{\theta}{c}\right)_{\#} p_\infty$ in distribution as $\tau \to 1$. As convergence with respect to $W_2$ implies weak convergence, the unique weak limit must be given as
\begin{equation}
     p_1^{\hat Z} = \left(\frac{\theta}{c}\right)_{\#} p_\infty.
\end{equation}
Letting $\tau_2 \to 1$ in \eqref{compressed-Z-Wasserstein} yields the endpoint estimate
\begin{equation}
 W_2(p_\tau^{\hat Z}, p_1^{\hat Z}) = \lim_{\tau_2 \to 1} W_2(p_\tau^{\hat Z}, p_{\tau_2}^{\hat Z}) \leq \lim_{\tau_2 \to 1} \tau_2 - \tau = 1 - \tau.
\end{equation}

\noindent (ii) 
Let $0 \le \tau_1 < \tau_2 \le 1$. Note that the case $\tau_2 = 1$ is already included in the estimate \eqref{Wasserstein-Bound}, using $\hat Y_\tau = Y_{-\frac{\ln(1-\tau)}{\theta}}$.
Hence, we are left to consider the case $\tau_2 < 1$. 
To this end, we adopt the notation $W_2(X,Y) \coloneqq W_2({\rm Law}(X), {\rm Law}(Y))$ for random variables $X,Y$. Then, it holds
\begin{equation}
    W_2\left(p_{\tau_1}^{\hat Y}, p_{\tau_2}^{\hat Y}\right) 
    = W_2\left( (1-\tau_1) Y_0 + \frac{c}{\theta} \hat Z_{\tau_1}, ~ (1-\tau_2) Y_0 + \frac{c}{\theta} \hat Z_{\tau_2}   \right).
\end{equation}
Since $Y_0$ and $\hat Z_{\tau_1}$ are independent, and $Y_0$ and $\hat Z_{\tau_2}$ are independent, it follows by the Wasserstein subadditivity under independence, see \cite[Corollary 1]{MariucciReiss2018}, that
\begin{equation}
    W_2\left(p_{\tau_1}^{\hat Y}, p_{\tau_2}^{\hat Y}\right) 
    \le W_2\big((1-\tau_1) Y_0, (1-\tau_2) Y_0\big) 
    + W_2\left(\frac{c}{\theta} \hat Z_{\tau_1}, \frac{c}{\theta} \hat Z_{\tau_2}\right).
\end{equation}
The first term satisfies
\begin{equation}
    W_2\big((1-\tau_1) Y_0, (1-\tau_2) Y_0\big) \le (\tau_2 - \tau_1)\mathbb{E}[|Y_0|^2]^{\frac{1}{2}},
\end{equation}
while the second term is given by the Wasserstein scaling property, see \cite[Lemma 2]{MariucciReiss2018}, as
\begin{equation}
    W_2\left(\frac{c}{\theta} \hat Z_{\tau_1}, \frac{c}{\theta} \hat Z_{\tau_2}\right) = \frac{c}{\theta} W_2\left(\hat Z_{\tau_1}, \hat Z_{\tau_2}\right).
\end{equation}
Now, the first part (i) delivers the claim.
\end{proof}

\subsection{Identification of Kac velocity fields}\label{app:odd-flux}

In general, velocity fields associated with a probability curve are not unique. Yet, imposing odd symmetry resolves this problem.

Let us first state an elementary symmetry property of the mode-wise Kac measures $\rho_{i,t}$.
\begin{lemma}\label{lem:odd-modewise}
    Let $(Z_t, I_t)_{t}$ be an extended, symmetric and inhomogeneous Kac process with speed $c_Z > 0$ and switching rate $a(t)$, constructed as in Proposition \ref{Z-Kac-is-MSP}. Let $\rho_{i, t}$, $i=1,2$, be the mode-wise subprobability measures \eqref{mode-wise-measure}. Then it holds the symmetry
    \begin{equation}
        \rho_{i, t}(A) = \rho_{j, t}(-A) \qquad \text{for all } A \in \mathcal{B}(\R), ~ i \neq j.
    \end{equation}
\end{lemma}

\begin{proof}   
    Consider the map
    \begin{equation}
        F: \{-1, +1\} \times D([0, \infty); \N) \to D([0, \infty); \R \times \{1, 2\}),\ (i, g) \mapsto \left( ci\int_0^t(-1)^{g(s)}\d s, \frac{1}{2}i(-1)^{g(t)} + \frac{3}{2}\right)_{t \geq 0},
    \end{equation}
    where $D([0, \infty); X)$ denotes the Skorohod space equipped with the $J_1$-topology.
    It is standard but technical to show that $F$ is measurable; we omit the details.
    By \eqref{inhom-Kac-process-Zt} and \eqref{switching-processes-from-poisson}, we can rewrite
    \begin{equation}
        (Z_t, I_t)_{t} = F(B_{\frac{1}{2}}, (N_t)_{t}).
    \end{equation}
    Since $B_{\frac{1}{2}} \overset{d}{=} - B_{\frac{1}{2}}$ and $F$ is measurable, we obtain
    \begin{equation}\label{flipped-distributional-equivalence}
        (Z_t, I_t)_{t} \overset{d}{=} F(-B_{\frac{1}{2}}, (N_t)_{t}) 
        = (-Z_t, 3- I_t)_{t},
    \end{equation}
    where the last equality holds pathwise by construction \eqref{inhom-Kac-process-Zt} and \eqref{switching-processes-from-poisson}.
    Now, for any $A \in \mathcal{B}(\R)$, it follows 
    \begin{equation}
        \rho_{i, t}(A) = \mathbb{P}(Z_t \in A, I_t = i) 
        = \mathbb{P}(Z_t \in -A, I_t = 3-i) = \rho_{j, t}(-A), \qquad i,j =1,2, ~ i \ne j.
    \end{equation}
\end{proof}

This leads to the odd symmetry of the Kac velocity field tied to the PDMP setting.
\begin{proposition}\label{prop:odd-velo}
    Let $(Z_t)_{t}$ be a symmetric and inhomogeneous Kac process with speed $c_Z > 0$ and switching rate $a(t)$, and let $p_t$ denote the law of $Z_t$.
    Then, the velocity field   
    \begin{equation}\label{Kac-expectation-velo}
        v_t(x) \coloneqq \mathbb{E}[ \dot Z_t \mid Z_t = x]
    \end{equation}
    is associated with $(Z_t)_t$, and $x\mapsto v_t(x)$ is an odd function (outside a $p_t$-null set).
\end{proposition}
\begin{proof}
    Proposition \ref{Z-Kac-is-MSP} and Proposition \ref{MSP-velo} immediately yield that $v_t$ is a velocity field associated with $(Z_t)_t$. Furthermore, equation \eqref{msp-velocity-formula} takes the form
    \begin{equation}
        v_t(x) = \sum_{i = 1}^2 (-1)^i c_Z \frac{\d \rho_{i, t}}{\d p_t}(x)
        = \frac{\d \rho_{2, t}}{\d p_t}(x) - \frac{\d \rho_{1, t}}{\d p_t}(x).
    \end{equation}
    By the uniqueness of the Radon-Nikodym derivative, Lemma \ref{lem:odd-modewise} implies $\frac{\d \rho_{2, t}}{\d p_t}(x) = \frac{\d \rho_{1, t}}{\d p_t}(-x)$ for $p_t$-almost every $x \in \R$, from which we infer that
    \begin{equation}
        v_t(x)
        = \frac{\d \rho_{2, t}}{\d p_t}(x) - \frac{\d \rho_{1, t}}{\d p_t}(x)
        = \frac{\d \rho_{1, t}}{\d p_t}(-x) - \frac{\d \rho_{2, t}}{\d p_t}(-x)
        = - v_t(-x)
    \end{equation}
     for $p_t$-almost every $x \in \R$. 
\end{proof}
Now, we can identify the following velocity fields associated with $(Z_t)_t$.
\begin{corollary}
    Let $(Z_t)_{t}$ be a symmetric and inhomogeneous Kac process with speed $c_Z=1$, and $p_t$ be the law of $Z_t$. Let $v_t^{\rm MSP}$ denote the velocity field \eqref{Kac-expectation-velo}, and 
    $v_t^{Z}$ denote the velocity field \eqref{Z-velocity}. Then, $v_t^{\rm MSP}$ and $v_t^{Z}$ coincide up to a $p_t$-null set.
\end{corollary}
\begin{proof}
    Fix $t > 0$. Now, both $v_t^{\rm MSP}$ and $v_t^{Z}$ satisfy the continuity equation \eqref{whole-CE} and are odd-symmetric by Prop.\ \ref{prop:odd-velo} and \eqref{Z-velocity}. Moreover, $p_t$ is even-symmetric and has an absolutely continuous part on $(-t,t)$ by \eqref{legendre-formula-total}. Altogether, for a.e.\ $t > 0$, we have
    \begin{equation}
       2 \, p_t(x) \, v_t^Z(x) = \int_{-x}^x \partial_x \big(p_t \, v_t^Z \big) \d x
       = - \int_{-x}^x \partial_t p_t \d x
       = \int_{-x}^x \partial_x \big(p_t \, v_t^{\rm MSP} \big) \d x
       = 2 \, p_t(x) \, v_t^{\rm MSP}(x),
    \end{equation}
    for ($p_t$-)almost all $|x| < t$, 
    which implies $v_t^Z(x) = v_t^{\rm MSP}(x)$ for ($p_t$-)a.e.\ $x \in (-t,t)$.
    Finally, the exact arguments as in the last part of the proof of Theorem \ref{cor:velo-formula} can be applied to $v_t^{\rm MSP}$, showing that $v_t^{\rm MSP}(\pm t) = \pm 1 = v_t^{Z}(\pm t)$, which finishes the proof.
\end{proof}

\section{Appendix: Numerical experiments}\label{sec:appendix-B}
\subsection{Metrics for discrete data}\label{sec:metric-discrete}

In Section \ref{sec:sharpness}, we test our OU-Kac flow matching model on different discrete datasets to ascertain its performance on very sharp datasets. Let $n \in \N$ denote the number of datapoints and $\mathcal{D} \coloneq \{x^{(i)} \in \R^{d} \mid i = 1, \ldots, n\}$ the dataset. We consider discrete distributions of the form
\begin{equation}\label{eq:discrete-measure}
    P_{\mathcal{D}} \coloneq \sum_{i = 1}^{n} p_{i} \delta_{\{x^{(i)}\}},
\end{equation}
where $p_{i} \geq 0,\ i = 1, \ldots, n$ and $\sum_{i = 1}^{n} p_{i} = 1$. We denote the learned distribution of the generative model by $\hat{P}_{\mathcal{D}}$. 

To describe the sharpness of our generative methods we define the relative distance (from now on distance) to $\mathcal{D}$ by 
\begin{equation*}
    d: \R^{d} \to [0, \infty),\ y \mapsto \frac{\min_{i} \|y - x^{(i)}\|}{\min_{i \neq j} \|x^{(i)} - x^{(j)}\|}.
\end{equation*}
For a sphere of radius $1$, we consider 
\begin{equation}
    d: \R^d \to [0, \infty),\ y \mapsto |\|y\| - 1|.
\end{equation}
For the rest of this section, let $\mathcal{S} = \{s^{(1)}, \ldots, s^{(N)}\} \subset \R^{d}$ be the set of generated samples. For every $i = 1, \ldots, n$, we define the \emph{recovered mass around mode} $x^{(i)} \in \mathcal{D}$ by 
\begin{equation*}
    R_{i}: [0, \infty) \to [0, 1],\ \varepsilon \mapsto \frac{1}{N} \sum_{k = 1}^{N} \mathds{1}_{\{\argmin_{j}\|s^{(k)} - x^{(j)}\| = i,\ d(s^{(k)}) \leq \varepsilon\}}
\end{equation*}
and the \emph{total recovered mass} by
\begin{equation*}
    R: [0, \infty) \to [0, 1],\ \varepsilon \mapsto \frac{1}{N} \sum_{k = 1}^{N} \mathds{1}_{\{d(s^{(k)}) \leq \varepsilon\}}.
\end{equation*}
A total recovery of $R(0.02) = 0.8$ means 80\% of generated samples have a distance of less than $0.02$ to the true data. Note that $R(\varepsilon) = \sum_{i = 1}^{n} R_{i}(\varepsilon)$ and for $S \sim \hat{P}_{\mathcal{D}}$ it holds
\begin{equation*}
    R(\varepsilon) = \hat{P}_{\mathcal{D}}(d(S) \leq \varepsilon).
\end{equation*}
The total recovered mass does not account for mode coverage. Indeed, if all $N$ generated samples have distance less than 0.02 to one fixed $x_{i} \in \mathcal{D}$, the total recovery would be 100\%, although the model has collapsed to one mode. To penalize such collapses, we consider the \emph{balanced recovery}
\begin{equation*}
    R_{\operatorname{balanced}}: [0, \infty) \to [0, 1],\ \varepsilon \mapsto \sum_{i = 1}^{n} \min \{p_{i}, R_{i}(\varepsilon)\},
\end{equation*}
where $p_{i}$ are the mode-probabilities from \eqref{eq:discrete-measure}. The balanced recovery of a very sharp, but mode collapsed model could only ever reach the true probability of the mode it collapsed to. Displaying $R(\varepsilon)$ or $R_{\operatorname{balanced}}(\varepsilon)$ only for some fixed $\varepsilon$ can be misleading. For example, setting $\varepsilon = 0.02$, a sample having relative distance 0.03 is considered equally bad as a sample having distance 1. Hence we further introduce the \emph{area under curve (AUC) sharpness} profile as the logarithmically weighted average of recovery profiles. Given lower and upper tolerance bounds $0 < a < b < \infty$, we let
\begin{equation}\label{sharpness_AUC}
    \operatorname{Sharp}_{\operatorname{AUC}}(a, b) \coloneq \frac{1}{\log(b) - \log(a)} \int_{a}^{b} \frac{R(\varepsilon)}{\varepsilon} \d{\varepsilon}.
\end{equation}
Equivalently, considering the weighting function
\begin{equation}
    w: [0, \infty) \to [0, 1],\ y \mapsto \begin{cases}1 & y \leq a\\ \frac{\log(b) - \log(y)}{\log(b) - \log(a)}, & y \in (a, b)\\ 0, & y \geq b\end{cases}
\end{equation}
we obtain 
\begin{equation}
    \operatorname{Sharp}_{\operatorname{AUC}}(a, b) = \frac{1}{N} \sum_{k = 1}^{N} w(d(s^{(k)})).
\end{equation}
By logarithmic weighting, improvements by a multiplicative constant are weighted equally. For example decreasing the distance from $2$ to $1$ is weighted the same as reducing it from $0.1$ to $0.05$ since $\log(2) - \log(1) = \log(0.1) - \log(0.05) = \log(2)$. This metric is again not aware of mode coverage, which is why we further consider the \emph{balanced AUC sharpness} profile given by
\begin{equation}\label{balanced_sharpness_AUC}
    \operatorname{Sharp}_{\operatorname{balanced, AUC}}(a, b) \coloneq \frac{1}{\log(b) - \log(a)} \int_{a}^{b} \frac{R_{\operatorname{balanced}}(\varepsilon)}{\varepsilon} \d{\varepsilon}.
\end{equation}

\subsection{Implementation details}\label{sec:implementation-details}
In this section, we provide further details on the numerical setup and experiments presented in Section \ref{sec:generative-modeling}. 
Our code is available online\footnote{https://github.com/auffenberg/finite-speed-ou-flows}.
\paragraph{OU-K model.}
Each OU-K model $(\tilde Y_t)_{t\in [0,1]}$ has parameters $a, c, \theta$ and $\mu$. By \eqref{invariant-OU-K}, the latent distribution has coordinate-wise variance
\begin{equation}
    \sigma_{{\rm OU-K}}^2 \coloneq \operatorname{Var}(\tilde Y_1) = \frac{c^2}{\theta(2 a + \theta)}.
\end{equation}
In all experiments, we consider the mean $\mu = 0$ and mean-reversion rate $\theta = 1$. 
Given a damping parameter $a$, we adjust the Kac speed $c$ to match the latent variance $\sigma_{{\rm OU-K}}^2$ to the OU baseline latent variance $\sigma_{{\rm OU}}^2$, which is specified for each experiment.

\paragraph{OU diffusion baseline.}
For the OU-diffusion baseline, we sample training times uniformly from $\operatorname{Unif}([10^{-5}, 1))$ and stop the reverse-time integration at network time $10^{-3}$. Throughout our experiments, we set $\theta=1$, and choose the OU volatility parameter $\sigma$ such that the OU terminal variance $\sigma_{\rm OU}^2 = \frac{\sigma^2}{2\theta}$ lies in $\{\frac{1}{100}, 1\}$.

\paragraph{Kac baseline.}
We use the Kac process from \cite[Section 7.2]{duong2026telegraphersgenerativemodelkac} with schedules $(1 - t, t^2)$ running up to the terminal time $T = 1$, that is 
\begin{equation}
    X_t = (1 - t)X_0 + K_{t^2}, \qquad t \in [0, 1],
\end{equation}
where the Kac process $(K_t)_{t \in [0, 1]}$ is independent of $X_0$. The variance at time $t \in [0, 1]$ is 
\begin{equation}
    \operatorname{Var}(K_t) = \frac{c^2}{a} \left(t - \frac{1 - e^{-2 at}}{2a}\right).
\end{equation}
Given a switching rate $a>0$, we match the Kac latent variance $\operatorname{Var}(K_1)$ to the FM baseline latent variance $\sigma_{\rm FM}^2$ by a suitable adjustment of $c$.

\paragraph{Toy datasets in Section \ref{sec:sharpness}.}
Although we eventually fix the dimension at $d=128$, we keep it arbitrary for now to investigate its effect on the numerical setup.
\begin{itemize}
    \item \textbf{Latent variance:} Here, we fix the componentwise latent variance of all compared models at $\sigma_\text{latent}^2 = \sigma_{\rm FM}^2 = \sigma_{\rm OU}^2 = 1$, independently from the dimension $d$.
    \item \textbf{Datasets:} 
    For $r>0$, we consider the $d$-dimensional rescaled toy datasets 
\begin{itemize}
    \item[(a)] Cross-polytope: $D^{(1)}(r) = \{\pm r e_1, \dots, \pm r e_d\}$, with the unit vectors $e_1, \dots, e_d \in \R^d$.
    \item[(b)] Rotated cross-polytope: $D^{(2)}(r) = U D^{(1)}(r)$, with some orthogonal matrix $U$.
    \item[(c)] Simplex: $D^{(3)}(r) = \{x_0, \ldots, x_d\}$, where 
    \begin{equation}
        (x_k)_j = r \sqrt{\frac{d + 1}{d}}\begin{cases}
            \frac{1}{\sqrt{j(j + 1)}}, & \text{for } j > k,\\ 
            -\frac{j}{\sqrt{j(j + 1)}}, & \text{for } j = k,\\ 
            0 & \text{for } j < k,
        \end{cases} \qquad \text{for } j = 1, \ldots, d.
    \end{equation}
    These points define the corners of a regular simplex centered around the origin with
    \begin{equation}
        \sum_{k = 0}^d x_k = 0, \qquad \|x_k\|_2^2 = r^2 \qquad \text{and} \qquad \langle x_k, x_\ell \rangle = - \frac{r^2}{d}
    \end{equation}
    for all $k, \ell = 0, \ldots, d$ with $k \neq \ell$.
    
    \item[(d)] Sphere of radius $1$: $D^{(4)}(1) = D_4 = \{x \in \R^d: \|x\|_2 = 1\}$.
\end{itemize}

For each dataset, we consider a uniform target distribution $X^{(i)} \sim \operatorname{Unif}(D^{(i)}(r))$, $i = 1,2,3,4$,
whose covariance matrix is given by
\begin{equation}
    \operatorname{Cov}(X^{(i)}) = \frac{r^2}{d}I_d, \quad i=1,2,3, \qquad \text{and }  \operatorname{Cov}(X^{(4)}) = \frac{1}{d} I_d.
\end{equation}
Hence, the coordinate-wise data variances are given as $\sigma_{\text{data}}^2 \coloneq \frac{r^2}{d}$ for $i=1,2,3$, while for $X^{(4)}$, it is $\sigma_{\text{data}}^2 = \frac{1}{d}$. 
Since we keep the component-wise latent variance fixed at $\sigma_\text{latent}^2 = 1$, we also fix the component-wise data variance $\sigma^2_\text{data}$ independently from the dimension $d$ for $i=1,2,3$. Here, we choose $\sigma_\text{data}^2 = \frac{1}{21}$ for all discrete toy datasets $D_1,D_2,D_3$,  resulting in the rescaling parameter $r = 2.47$ in dimension $d=128$.
    
For the spherical dataset $D_4$, we retain the unit radius instead of applying some variance normalization, allowing us to better distinguish the behaviour of the compared models for this dataset.

\item\textbf{Metrics:} The main metric we consider is the balanced AUC sharpness  \eqref{balanced_sharpness_AUC} with $a = 0.005$ and $b = 0.1$, applied to the data sets $D_i, \, i=1,2,3$. 
Since the construction \eqref{balanced_sharpness_AUC} only makes sense for a discrete set with a finite number of elements, we apply the AUC sharpness \eqref{sharpness_AUC} for the sphere $D_4$.
To compute the integrals in \eqref{sharpness_AUC} and \eqref{balanced_sharpness_AUC}, we use the trapezoidal rule for $u = \log(\varepsilon)$, with 200 equidistantly spaced points in $u$.
Then, every 40,000 training steps, we generate 10,000 samples using the EMA model (see below) to compute the (balanced) AUC sharpness reported in Figure \ref{fig:four-sets-d128}.

\item \textbf{Neural Model:} We train a fully connected multilayer perceptron with 5 hidden layers of width 1472, SiLU activations and a final linear output layer of the dimension of the dataset for 400k iterations. The input is the dimension of the dataset with appended time. 

\item \textbf{Learning Hyperparameters:} We use the AdamW optimizer with Adam $\beta_1 = 0.9$ and Adam $\beta_2 = 0.99$, a learning rate of $10^{-4}$, a batch size of 2048, weight decay of $10^{-4}$, gradient clipping of 5 and an EMA of 0.9999. 

\item \textbf{Integration:} We use the classic Runge-Kutta (RK4) method with 100 RK4 steps for solving the ODE \eqref{eq:flow_ode}. Each RK4 step requires 4 neural function evaluations (NFE). In this experiment, we opted for RK4 instead of Euler integration to reduce discretization errors.
\end{itemize}

\paragraph{CIFAR-10 experiment in Section \ref{sec:straight-trajectories}.}
\begin{itemize}
    \item \textbf{Latent variance:} Here, we fix the componentwise latent variance of the FM and Kac baselines to $\sigma_{\rm FM}^2 = 1$, and of the OU and OU-K models to $\sigma_{\rm OU}^2 = \frac{1}{100}$. We have noticed that small latent variances are beneficial for the OU variants.
    \item \textbf{Schedule:} 
    For the OU-K model, we use the cosine schedule exactly as given in \eqref{cosine-schedule}, without any time truncation. For the OU baseline, we additionally apply beta-clipping at $0.999$ as in \cite{nichol2021improved}. Let $N = 1000$ and set 
    \begin{equation}
        \beta_i = \min \left\{1 - e^{-g\left(\frac{i + 1}{N}\right) + g\left(\frac{i}{N}\right)} , 0.999 \right\}, \qquad i = 0, \ldots, N - 1.
    \end{equation}
    Then we consider the schedule
    \begin{equation}
        \tilde g(t) = - \left(\sum_{j = 0}^{i - 1} \log(1 - \beta_j) + (N t - i)\log(1 - \beta_i)\right), \qquad t \in \left[\frac{i}{N}, \frac{i + 1}{N}\right].
    \end{equation}
        
    \item  \textbf{Model:} We train a time-conditioned U-Net architecture from \cite{nichol2021improved} for 400k iterations. The U-Net has 128 base channels, channel multipliers (1, 2, 2, 2), three residual blocks per level, 4 attention heads with 64 channels per head at spatial resolutions $16 \times 16$ and in the $4 \times 4$ bottleneck, a dropout of 0.1, scale-shift normalization and convolutional resampling.

    \item \textbf{Learning hyperparameters:} We use a learning rate of $2 \times 10^{-4}$, the Adam optimizer with Adam $\beta_1 = 0.9$ and $\beta_2 = 0.999$, a batch size of 128, no weight-decay, gradient clipping at 1, an EMA decay of 0.9999, a 5k step warmup where the learning rate increases linearly, BF16 autocast format and the fixed training seed 0. We extend the 50k dataset with a resolution of $32 \times 32$ by random horizontal flipping with probability $\frac{1}{2}$. The images are normalized to $[-1, 1]$ under the transformation $x \mapsto 2x - 1$.
    \item \textbf{Evaluation:} After 400k iterations we generate 50k images and compute the Fr{\`e}chet inception distance (FID) against the CIFAR-10 training set. We use the standard $32 \times 32$ image resolution and convert generated images to 8-bit RGB values. The ODE \eqref{eq:flow_ode} is solved using fixed-step Euler 
    integration with uniformly spaced time steps, where we showcase the results for 100 and 200 Euler steps. 
\end{itemize}

\newpage
\subsection{CIFAR-10 generated images}\label{subsec:cifar-figures}
For Section \ref{sec:straight-trajectories}, we showcase the generated CIFAR-10 images of our OU-K model and the baselines.
\begin{figure}[H] 
  \centering
  \begin{minipage}[b]{0.45\textwidth}
    \centering
    \includegraphics[width=\linewidth]{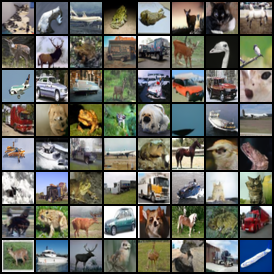}
    \par\smallskip
   { \( \rm FM\)}
  \end{minipage}
  \hfill
  \begin{minipage}[b]{0.45\textwidth}
    \centering
    \includegraphics[width=\linewidth]{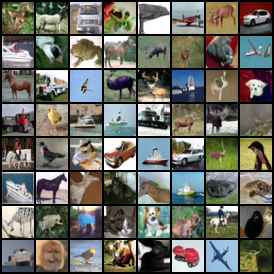}
    \par\smallskip
    {Kac}
  \end{minipage}
  \hfill \\[7ex]
  \begin{minipage}[b]{0.45\textwidth}
    \centering
    \includegraphics[width=\linewidth]{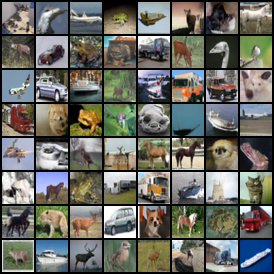}
    \par\smallskip
    {\( \rm OU\)}
  \end{minipage}
  \hfill
  \begin{minipage}[b]{0.45\textwidth}
    \centering
    \includegraphics[width=\linewidth]{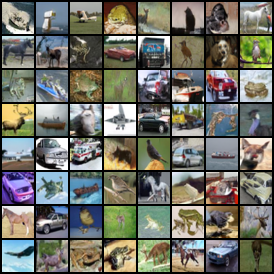}
    \par\smallskip
    {OU-K (ours)}
  \end{minipage}

  \caption{CIFAR-10 generation results of the FM, OU and Kac baselines and the OU-K model. }
  \label{fig:cifar-1}
\end{figure}


\end{document}